\documentclass{article}
\usepackage{amssymb,amsmath,bm,hyperref,rotating}
\usepackage{amssymb,amsmath,bm,geometry,graphics,graphicx,url,color,inputenc}
\def\no{\noindent}
\def\pmatrix{\left(\begin{array}}
\def\endpmatrix{\end{array}\right)}

\def\red{\color{black}}
\def\blue{\color{black}}

\def\RR{\mathbb{R}}

\def\B{{\cal B}}

\def\F{{\cal F}}

\def\I{{\cal I}}

\def\P{{\cal P}}

\def\dd{\mathrm{d}}

\newtheorem{theo}{Theorem}
\newtheorem{lem}{Lemma}
\newtheorem{cor}{Corollary}
\newtheorem{rem}{Remark}
\newtheorem{defi}{Definition}
\def\proof{\noindent\underline{Proof}\quad}
\def\QED{\mbox{~$\Box{~}$}}

\def\bfb{{\bm{b}}}
\def\bfc{{\bm{c}}}

\def\bfe{{\bm{e}}}

\def\bfzero{{\bm{0}}}

\def\bfgamma{{\bm{\gamma}}}
\def\bfdelta{{\bm{\delta}}}
\def\bfeta{{\bm{\eta}}}

\def\bfphi{{\bm{\phi}}}
\def\bfpi{{\bm{\pi}}}

\begin{document}

\title{Arbitrarily high-order energy-conserving methods for Poisson problems\,\thanks{Cite this article as: Amodio, P., Brugnano, L. \& Iavernaro, F. Arbitrarily high-order energy-conserving methods for Poisson problems. Numer Algor (2022) \url{https://doi.org/10.1007/s11075-022-01285-z}}}

\author{Pierluigi Amodio\footnote{Dipartimento di Matematica, Universit\`a di Bari, Italy. \qquad\, \url{{pierluigi.amodio,felice.iavernaro}@uniba.it}} \and Luigi Brugnano\footnote{Dipartimento di Matematica e Informatica ``U.\,Dini'', Universit\`a di Firenze, Italy. \quad\url{luigi.brugnano@unifi.it}}  \and  Felice Iavernaro$^*$}


\maketitle

\begin{abstract}  In this paper we are concerned with energy-conserving methods for Poisson problems, which are effectively solved by defining a  suitable generalization of HBVMs, a class of energy-conserving methods for Hamiltonian problems. The actual implementation of the methods is fully discussed, with a particular emphasis on the conservation of Casimirs. Some numerical tests are reported, in order to assess the theoretical findings.

\medskip
\no{\bf Keywords:~} Poisson problems, Casimir function, Line Integral Methods, Hamiltonian Boundary Value Methods, HBVMs.

\medskip
\no{\bf MSC:~} 65L05, 65P10.

\end{abstract}

\section{Introduction}\label{intro} A Poisson problem is in the form
\begin{equation}\label{poisson}
\dot y = B(y)\nabla H(y) =: f(y), \quad t>0,\qquad y(0)=y_0\in\RR^m, \qquad B(y)^\top = -B(y),
\end{equation}
where $H(y)$ is a scalar function, usually called the {\em Hamiltonian}. For sake of simplicity, hereafter both $H(y)$ and $B(y)$ will be assumed suitably regular. From the skew-symmetry of $B(y)$ one easily deduces that $H(y)$ is a  constant of motion, since, along the solution of (\ref{poisson}),
$$\frac{\dd}{\dd t} H(y) = \nabla H(y)^\top \dot y = \nabla H(y)^\top B(y)\nabla H(y) = 0,$$
due to the skew-symmetry of $B(y)$. Possible additional invariants of (\ref{poisson}) are its {\em Casimirs}, namely scalar functions $C(y)$ for which 
\begin{equation}\label{casimir}
\nabla C(y)^\top B(y)={\red (0,\dots,0)\in\RR^{1\times m}}
\end{equation}
holds true for all $y$. When matrix $B(y)$ is constant, as in the case of Hamiltonian problems,
\begin{equation}\label{Ham}
 \dot y = J\nabla H(y), \qquad t>0, \qquad y(0)=y_0, \qquad J^\top=-J,
 \end{equation}
energy conservation can be obtained by solving problem (\ref{Ham}) via HBVMs, a class of energy-conserving Runge-Kutta methods for Hamiltonian problems (see, e.g., \cite{BIT2009,BIT2010,BIT2012,BIT2015} and the monograph \cite{LIMbook2016}, see also the recent review paper \cite{BI2018}). Nevertheless, in the case where the problem is not Hamiltonian, HBVMs are no more energy-conserving. This motivates the present paper, where an energy-conserving variant of HBVMs for Poisson problems is derived and analyzed.

The numerical solution of Poisson problems has been tackled by following many different approaches (see, e.g., \cite[Chapter\,VII]{GNI2006} and references therein). More recently, it has been considered in \cite{CH2011}, where an extension of the AVF method \cite{QML2008} is proposed, and in \cite{BCMR2012,BGI2018}, where a line integral approach has been used instead. Functionally fitted methods have been proposed in {\red \cite{M2015,WW2018,WW2021,MHW2022}.} In this paper we further pursue the line integral approach to the problem, which will provide an energy-conserving variant of HBVMs for solving (\ref{poisson}).

With this premise, the structure of the paper is the following: in Section~\ref{new} we describe the new framework, in which the methods will be derived; in Section~\ref{discr} we provide the final shape of the method, while in Section~\ref{discrete} its actual implementation is studied; in Section~\ref{num} we present a few numerical tests confirming the theoretical findings; at last, in Section~\ref{fine} we give some concluding remarks.

\section{The new framework}\label{new}
As anticipated above, the framework that we shall use to derive and analyze the methods is that of the so called {\em line integral methods}, namely methods where the conservation properties are derived by the vanishing of a corresponding {\em line integral} \cite{LIMbook2016,BI2018}. Such methods have been largely investigated in the case of Hamiltonian problems, providing their major instance given by Hamiltonian Boundary Value Methods (HBVMs). The analysis will strictly follow that in \cite{BIT2012} and \cite{BMR2019}. To begin with, let us consider problem (\ref{poisson}) on the interval $[0,h]$,
\begin{equation}\label{poisson1c}
\dot y(ch) =  B(y(ch)) \nabla H(y(ch)),\qquad c\in[0,1], \qquad y(0) = y_0.
\end{equation}
In fact, since we shall speak about a one-step method, it suffices to analyse its first step of application, with $h$ the time-step. Next, let us consider the  orthonormal Legendre polynomial basis $\{P_j\}_{j\ge0}$ on the interval $[0,1]$,
\begin{equation}\label{orto}
\deg P_j = j, \qquad \int_0^1 P_i(x)P_j(x)\dd x = \delta_{ij}, \qquad \forall i,j=0,1,\dots,
\end{equation}
with $\delta_{ij}$ the Kronecker symbol, and the following expansions for the functions at the right-hand side in (\ref{poisson1c}):

\begin{eqnarray}\nonumber
\nabla H(y(ch)) = \sum_{j\ge0} P_j(c) \gamma_j(y), &&P_j(c)B(y(ch))  = \sum_{i\ge0} P_i(c) \rho_{ij}(y), \qquad~  c\in[0,1],\\[-3mm] \label{gammaro}
\\[-3mm]\nonumber
\gamma_j(y)     = \int_0^1 P_j(\tau)\nabla H(y(\tau h))\dd\tau, &&
\rho_{ij}(y)         = \int_0^1 P_i(\tau)P_j(\tau)B(y(\tau h))\dd\tau, \quad i,j=0,1,\dots. 
\end{eqnarray}

The following properties hold true.

\begin{lem}\label{Ohj} Assume $\psi:[0,h]\rightarrow V$, with $V$ a vector space, admit a Taylor expansion at 0. Then, for all $j=0,1,\dots$ :
$$\int_0^1 P_j(c)c^i\psi(ch)\dd c=O(h^{j-i}), \qquad i=0,\dots,j.$$\end{lem}
\proof 
By the hypotheses on $\psi$, one has: $$c^i\psi(ch)=\sum_{r\ge0} \frac{\psi^{(r)}(0)}{r!} h^r c^{r+i}.$$ 
Consequently, for all $i=0,\dots,j$, by virtue of (\ref{orto}) it follows that:
$$
\int_0^1 P_j(c)c^i\psi(ch)\dd c= \sum_{r\ge0} \frac{\psi^{(r)}(0)}{r!} h^r \int_0^1P_j(c)c^{r+i}\dd c
=\sum_{r\ge j-i} \frac{\psi^{(r)}(0)}{r!} h^r \int_0^1P_j(c)c^{r+i}\dd c = O(h^{j-i}).\QED 
$$               

\begin{cor}\label{grhj} With reference to (\ref{gammaro}), for any suitably regular path  $\sigma:[0,h]\rightarrow\RR^m$  one has: 
\begin{equation}\label{gammaj}
\gamma_j(\sigma)=O(h^j), \qquad \rho_{ij}(\sigma)=O(h^{|i-j|}).\qquad \forall i,j=0,1,\dots.
\end{equation}
\end{cor}
\proof Immediate from Lemma~\ref{Ohj}, by taking into account  (\ref{gammaro}).\,\QED\bigskip

We also state, without proof, the following straightforward property, deriving from the skew-symmetry of $B$.

\begin{lem}\label{lem_roij}
With reference to (\ref{gammaro}), for any path  $\sigma:[0,h]\rightarrow\RR^m$  one has: 
\begin{equation}\label{roij}
\rho_{ij}(\sigma) = \rho_{ji}(\sigma) = -\rho_{ij}(\sigma)^\top, \qquad \forall i,j=0,1,\dots.
\end{equation}
\end{lem} 

Taking into account (\ref{gammaro}),  the right-hand side in (\ref{poisson1c}) can be rewritten as:

\begin{equation}\label{infiy1}
\dot y(ch) = B(y(ch))\nabla H(y(ch)) = \sum_{j\ge0} P_j(c) B(y(ch)) \gamma_j(y)
=  \sum_{i,j\ge 0} P_i(c)\rho_{ij}(y) \gamma_j(y), \qquad c\in[0,1],
\end{equation}
from which one obtains that the solution of (\ref{poisson1c}) can be formally written as:
\begin{equation}\label{infiy}
y(ch) = y_0 + h\sum_{i,j\ge0} \int_0^cP_i(x)\dd x\rho_{ij}(y) \gamma_j(y), \qquad c\in[0,1].
\end{equation}
In particular, by considering (\ref{orto}) and that $P_0(c)\equiv 1$, from which $\int_0^1P_i(x)\dd x=\delta_{i0}$, one has:
\begin{equation}\label{infiyh}
y(h) = y_0 + h\sum_{j\ge 0}\rho_{0j}(y) \gamma_j(y) 
\equiv y_0 + h\sum_{j\ge0} \int_0^1 P_j(c)B(y(ch))\dd c\int_0^1 P_j(c) \nabla H(y(ch))\dd c. ~~
\end{equation}
In order to obtain a polynomial approximation of degree $s$ to $y$, it suffices to truncate the two infinite series in (\ref{infiy1}) after $s$ terms:
\begin{equation}\label{sigma1}
\dot \sigma(ch) =   \sum_{i,j=0}^{s-1} P_i(c)\rho_{ij}(\sigma) \gamma_j(\sigma), \qquad c\in[0,1],
\end{equation}
with $\rho_{ij}(\sigma)$ and $\gamma_j(\sigma)$ defined according to (\ref{gammaro}) by formally replacing $y$ by $\sigma$.
Consequently, (\ref{infiy}) becomes
\begin{equation}\label{sigma}
\sigma(ch) = y_0 + h\sum_{i,j=0}^{s-1} \int_0^cP_i(x)\dd x\rho_{ij}(\sigma) \gamma_j(\sigma), \qquad c\in[0,1],
\end{equation}
providing the approximation
\begin{equation}\label{sigmah}
y_1:=\sigma(h) = y_0 + h\sum_{j=0}^{s-1}\rho_{0j}(\sigma) \gamma_j(\sigma)
 \equiv y_0 + h\sum_{j=0}^{s-1} \int_0^1 P_j(c)B(\sigma(ch))\dd c\int_0^1 P_j(c) \nabla H(\sigma(ch))\dd c,
\end{equation}
in place of (\ref{infiyh}).

\subsection{Interpretation of $\sigma$}\label{intersigma} We now provide an interesting interpretation of the polynomial approximation $\sigma$. For this purpose, let us rewrite (\ref{infiy1}), by taking into account (\ref{gammaro}), as follows:
\begin{eqnarray*}
\dot y(ch) &=& \sum_{i,j\ge0} P_i(c) \rho_{ij}(y)\gamma_j(y)\\
&=&\sum_{i\ge0} P_i(c) \sum_{j\ge0} \int_0^1 P_i(\tau) B(y(\tau h))P_j(\tau)\dd\tau\int_0^1 P_j(\tau_1)\nabla H(y(\tau_1h))\dd\tau_1\\
&=&\sum_{i\ge0} P_i(c)  \int_0^1 P_i(\tau) B(y(\tau h))\left(\underbrace{\sum_{j\ge0}P_j(\tau)\int_0^1 P_j(\tau_1)\nabla H(y(\tau_1h))\dd\tau_1}_{=\nabla H(y(\tau))}\right)\dd\tau\\
&\equiv& B(y(ch))\nabla H(y(ch)),
\end{eqnarray*}
as is expected. In a similar way, we can rewrite (\ref{sigma1}) as:
\begin{eqnarray*}
\dot \sigma(ch) &=& \sum_{i,j=0}^{s-1} P_i(c) \rho_{ij}(\sigma)\gamma_j(\sigma)\\
&=&\sum_{i=0}^{s-1} P_i(c) \sum_{j=0}^{s-1} \int_0^1 P_i(\tau) B(\sigma(\tau h))P_j(\tau)\dd\tau\int_0^1 P_j(\tau_1)\nabla H(\sigma(\tau_1h))\dd\tau_1\\
&=&\sum_{i=0}^{s-1} P_i(c)  \int_0^1 P_i(\tau) B(\sigma(\tau h))\left(\sum_{j=0}^{s-1}P_j(\tau)\int_0^1 P_j(\tau_1)\nabla H(\sigma(\tau_1h))\dd\tau_1\right)\dd\tau\\ &=:& \sum_{i=0}^{s-1} P_i(c)  \int_0^1 P_i(\tau) B(\sigma(\tau h))\left[\nabla H(\sigma(\tau h))\right]_s\dd\tau\\
&\equiv& \left[ B(\sigma(\tau h))\left[\nabla H(\sigma(\tau h))\right]_s\right]_s,
\end{eqnarray*}
having denoted by $\left[\cdot\right]_s$ the best approximation in  {\red $\Pi_{s-1}$ (i.e., $[\cdot]_s$ is the best polynomial approximation of degree $s-1$)} of the function in argument. This fact provides a noticeable interpretation of the polynomial approximation $\sigma$, which is the solution of the initial value problem
\begin{equation}\label{odeivp_s}
\dot\sigma(ch) = \left[B(\sigma(ch))\left[\nabla H(\sigma(ch))\right]_s\right]_s, \qquad c\in[0,1], \qquad \sigma(0)=y_0,
\end{equation}
equivalent to (\ref{sigma1}). Thus, the vector field of (\ref{odeivp_s}) is defined by a double projection procedure onto the finite dimensional vector space {\red $\Pi_{s-1}$} which involves, in turn, the vector fields $\nabla H(\sigma(ch))$ and 
$B(\sigma(ch))\left[\nabla H(\sigma(ch))\right]_s$, respectively.


\subsection{Analysis}\label{analsigma} We now analyze the method (\ref{sigma1})--(\ref{sigmah}). The following result then holds true, stating that the method is energy-conserving.

\begin{theo}\label{Hcons} $H(y_1)=H(y_0)$.
\end{theo}
\proof In fact, by virtue of (\ref{poisson}), (\ref{gammaro}), and (\ref{sigma1})--(\ref{sigmah}) one has, by using the standard line integral argument:
\begin{eqnarray*}
\lefteqn{H(y_1)-H(y_0)~=~ H(\sigma(h))-H(\sigma(0)) ~=~ \int_0^h \nabla H(\sigma(t))^\top\dot\sigma(t)\dd t}\\ 
&=& h\int_0^1 \nabla H(\sigma(ch))^\top\dot\sigma(ch)\dd c~=~ h\int_0^1 \nabla H(\sigma(ch))^\top \sum_{i,j=0}^{s-1} P_i(c)\rho_{ij}(\sigma) \gamma_j(\sigma)\dd c\\
 &=&h \sum_{i,j=0}^{s-1} \left[\int_0^1 P_i(c)\nabla H(\sigma(ch))\dd c\right]^\top\rho_{ij}(\sigma) \gamma_j(\sigma)
 ~=~h \sum_{i,j=0}^{s-1} \gamma_i(\sigma)^\top\rho_{ij}(\sigma) \gamma_j(\sigma)
~=~0,
\end{eqnarray*}
where the last equality follows from (\ref{roij}).\,\QED\bigskip

Concerning the accuracy of the new approximation, the following result holds true.

\begin{theo}\label{y1yh} Let $y_1$ be defined according to (\ref{sigma1})--(\ref{sigmah}). Then, $y_1-y(h)=O(h^{2s+1})$.\footnote{I.e., the approximation procedure has order of convergence $2s$.}\end{theo}

 \proof Let $y(t)\equiv y(t,\xi,\eta)$ denote the solution of the initial value problem (see (\ref{poisson}))
 \begin{equation}\label{poissonxi}
 \dot y = B(y) \nabla H(y) =:F(y), \qquad t\ge\xi, \qquad y(\xi)=\eta, 
 \end{equation}
Moreover, let us denote
{\red
$$\Phi(t,\xi,\eta) = \frac{\partial}{\partial \eta}y(t,\xi,\eta),\qquad t\ge \xi,$$
also recalling that 
$$\qquad\frac{\partial}{\partial \xi}y(t,\xi,\eta) = - \Phi(t,\xi,\eta)F(\eta).$$
 Then, by taking into account Lemma~\ref{Ohj} and Corollary~\ref{grhj}, and setting  
 $$\Psi_i(\sigma) = \int_0^1 P_i(c)\Phi(h,ch,\sigma(ch))\dd c =O(h^i), \qquad i=0,1,\dots,$$ }
 one has:
 \begin{eqnarray*}
 \lefteqn{
 y_1-y(h) ~=~\sigma(h)-y(h) = y(h,h,\sigma(h))-y(h,0,\sigma(0)) = \int_0^h \frac{\dd}{\dd t} y(h,t,\sigma(t))\dd t}\\
 &=&\int_0^h \left.\left[\frac{\partial}{\partial \xi} y(h,\xi,\sigma(t))\right|_{\xi=t} + \left.\frac{\partial}{\partial \eta}y(h,t,\eta)\right|_{\eta=\sigma(t)}
 \dot\sigma(t)\right] \dd t\\
 &=& \int_0^h \left[-\Phi(h,t,{\red \sigma(t)})F(\sigma(t))+\Phi(h,t,{\red \sigma(t)})\dot\sigma(t)\right]\dd t  \\
 &=& -h\int_0^1 \Phi(h,ch,{\red \sigma(ch)})\left[ F(\sigma(ch))-\dot\sigma(ch)\right]\dd c\\
 &=&-h\int_0^1 \Phi(h,ch,{\red \sigma(ch)})\left[ B(\sigma(ch))\nabla H(\sigma(ch)) - \sum_{i,j=0}^{s-1} P_i(c)\rho_{ij}(\sigma)\gamma_j(\sigma)\right]\dd c\\
 &=&-h\int_0^1 \Phi(h,ch,{\red \sigma(ch)})\left[ \sum_{i,j\ge 0} P_i(c)\rho_{ij}(\sigma)\gamma_j(\sigma) - \sum_{i,j=0}^{s-1} P_i(c)\rho_{ij}(\sigma)\gamma_j(\sigma)\right]\dd c\\
 &=& -h\left[\sum_{i,j\ge0} \Psi_i(\sigma)\rho_{ij}(\sigma)\gamma_j(\sigma) -\sum_{i,j=0}^{s-1} \Psi_i(\sigma)\rho_{ij}(\sigma)\gamma_j(\sigma) \right]\\
 &=&-h\left[ \sum_{i=0}^{s-1}\sum_{j\ge s}  \underbrace{\Psi_i(\sigma)\rho_{ij}(\sigma)}_{=O(h^j)}\gamma_j(\sigma) +  \sum_{i\ge s}\sum_{j=0}^{s-1}  \Psi_i(\sigma)\underbrace{\rho_{ij}(\sigma)\gamma_j(\sigma)}_{=O(h^i)}
 +  \sum_{i,j\ge s}  \Psi_i(\sigma)\rho_{ij}(\sigma)\gamma_j(\sigma)
 \right]\\[2mm] &=&O(h^{2s+1}).\,\QED
 \end{eqnarray*}\bigskip
 
 At last, we observe that, since the procedure (\ref{sigma1})--(\ref{sigmah}) is equivalent to defining the path $\sigma$ that joins $\sigma(0)=y_0$ to  $\sigma(h) = y_1$, then the same procedure, when started at $y_0$ and going forward provides $y_1$ and, when started from $y_1$ and going backward, brings back to $y_0$. In other words, the following result holds true.
 
 \begin{theo}\label{simme} The procedure (\ref{sigma1})--(\ref{sigmah}) is symmetric.\end{theo}  
\proof
This result comes as an easy consequence of Theorem \ref{simme1}, where the analogous property for the full discretized method is shown. \,\QED
 
\begin{rem}\label{HBVM1} We conclude this section emphasizing that, when problem (\ref{poisson}) is Hamiltonian, i.e., in the form (\ref{Ham}),
 then matrix $B(y)\equiv J$ is constant, and therefore (see (\ref{gammaro})) $\rho_{ij}(\sigma) = {\red \delta_{ij}}J$.
 Consequently, equation (\ref{sigma}) becomes
 $$\sigma(ch) = y_0+h\sum_{j=0}^{s-1} \int_0^c P_j(x)\dd x \int_0^1 P_j(\tau)J\nabla H(\sigma(\tau h))\dd\tau, \qquad c\in[0,1],$$
 which (see \cite{BIT2010,BIT2012-1,BI2018}) is the so called ``\,{\em master functional equation}'' defining the class of energy-conserving methods named {\em Hamiltonian Boundary Value Methods (HBVMs)}. Consequently, when the problem is Hamiltonian, then the procedure (\ref{sigma1})--(\ref{sigmah}) reduces to the {\em HBVM$(\infty,s)$} method in \cite{BIT2010}.
 \end{rem}
 
 \subsection{Conservation of Casimirs}\label{cassec} In this section, we study the required modifications, in order to conserve Casimirs, i.e., functions satisfying (\ref{casimir}). For sake of simplicity, we shall consider the simpler case of one Casimir, but multiple ones can be handled by slightly adapting the arguments, as is sketched at the end of the section. To begin with, for the original problem (\ref{poisson1c}), and its equivalent formulation (\ref{infiy1}), one has:
\begin{eqnarray}\nonumber
0 &=& C(y(h))-C(y_0) = \int_0^h \nabla C(y(t))^\top \dot y(t)\dd t = h\int_0^1 \nabla C(y(ch))^\top\dot y(ch)\dd c \\ \nonumber
&=&  h\int_0^1 \nabla C(y(ch))^\top B(y(ch))\nabla H(y(ch))\dd c =h\int_0^1 \nabla C(y(ch))^\top\sum_{i,j\ge0} P_i(c) \rho_{ij}(y)\gamma_j(y) \\
&=& h\sum_{i,j\ge0}\left[ \int_0^1 P_i(c)\nabla C(y(ch))\dd c\right]^\top
   \rho_{ij}(y)\gamma_j(y) ~=:~ h\sum_{i,j\ge0} \pi_i(y)^\top\rho_{ij}(y)\gamma_j(y),\label{casy}
\end{eqnarray}
having set
\begin{equation}\label{piy}
\pi_i(y) = \int_0^1 P_i(c)\nabla C(y(ch))\dd c = O(h^i),
\end{equation}
the $i$-th Fourier coefficient of the gradient of the Casimir, with the last equality following from Lemma~\ref{Ohj}. Clearly, again from (\ref{casimir}), one derives, by taking into account (\ref{sigma1}):
\begin{eqnarray}\nonumber
\lefteqn{C(y_1)-C(y_0) ~=~C(\sigma(h))-C(\sigma(0)) ~=~\int_0^h \nabla C(\sigma(t))^\top \dot \sigma(t)\dd t ~=~ h\int_0^1 \nabla C(\sigma(ch))^\top\dot \sigma(ch)\dd c}\\ \nonumber
&=& h\int_0^1 \nabla C(\sigma(ch))^\top\left[\dot \sigma(ch)-B(\sigma(ch))\nabla H(\sigma(ch))\right]\dd c\\ \nonumber
&&+~\overbrace{h\int_0^1 \nabla C(\sigma(ch))^\top B(\sigma(ch))\nabla H(\sigma(ch))\dd c}^{=0}\\ \nonumber
&=&-h\int_0^1  \nabla C(\sigma(ch))^\top\left[ \sum_{i,j\ge0} P_i(c) \rho_{ij}(\sigma)\gamma_j(\sigma) - \sum_{i,j=0}^{s-1} P_i(c) \rho_{ij}(\sigma)\gamma_j(\sigma)\right]\dd c\\ \nonumber
&=& -h\left[ \sum_{i,j\ge s}\pi_i(\sigma)^\top\rho_{ij}(\sigma)\gamma_j(\sigma) + \sum_{i=0}^{s-1}\sum_{j\ge s}\underbrace{\pi_i(\sigma)^\top\rho_{ij}(\sigma)}_{=O(h^j)}\gamma_j(\sigma)+ \sum_{j=0}^{s-1}\sum_{i\ge s}\pi_i(\sigma)^\top\underbrace{\rho_{ij}(\sigma)\gamma_j(\sigma)}_{=O(h^i)}\right]
\\&=&O(h^{2s+1}).\label{casig}
\end{eqnarray}

In order to recover the conservation of Casimirs, we shall use a strategy akin to that used in \cite{BS2014} for HBVMs (see also \cite{BI2012}), i.e., suitably perturbing some of its coefficients. In more details, let us consider the following modified polynomial in place of (\ref{sigma1}):\footnote{Here, we take into account that $P_0(c)\equiv1$.}
\begin{equation}\label{sigma1alfa}
\dot \sigma_\alpha(ch) =   \sum_{i,j=0}^{s-1} P_i(c)\rho_{ij}(\sigma_\alpha) \gamma_j(\sigma_\alpha) - \alpha \tilde{B}\gamma_0(\sigma_\alpha), \qquad c\in[0,1], \qquad \sigma_\alpha(0)=y_0, 
\end{equation}
with $\tilde{B}^\top=-\tilde{B}\ne O$ an arbitrary skew-symmetric matrix. As is usual, the new approximation will be $y_1:=\sigma_\alpha(h)$. In other words, we have considered the following perturbed coefficient:
$$\rho_{00}(\sigma_\alpha) - \alpha \tilde{B},$$
in place of $\rho_{00}(\sigma)$ in (\ref{sigma1}). The following result holds true.

\begin{theo}\label{thalfa}
Assume that $\pi_0(\sigma_\alpha)^\top 
\tilde{B}\gamma_0(\sigma_\alpha)\ne 0$. Then the Casimir $C(y)$ is conserved, provided that
\begin{equation}\label{alfa} 
\alpha = \frac{\sum_{i,j=0}^{s-1} \pi_i(\sigma_\alpha)^\top \rho_{ij}(\sigma_\alpha)\gamma_j(\sigma_\alpha)}{\pi_0(\sigma_\alpha)^\top \tilde{B}\gamma_0(\sigma_\alpha)}.
\end{equation}
Moreover, $\alpha = O(h^{2s})$.
\end{theo}
\proof In fact, by repeating similar steps as in (\ref{casig}), and replacing $\sigma$ by $\sigma_\alpha$, as defined in (\ref{sigma1alfa}), one obtains:
\begin{eqnarray*} 
\lefteqn{C(y_1)-C(y_0) ~=~C(\sigma_\alpha(h))-C(\sigma_\alpha(0))~=~h\int_0^1 \nabla C(\sigma_\alpha(ch))^\top\dot \sigma_\alpha(ch)\dd c}\\
&=&h\int_0^1  \nabla C(\sigma_\alpha(ch))^\top\left[\sum_{i,j=0}^{s-1} P_i(c) \rho_{ij}(\sigma_\alpha)\gamma_j(\sigma_\alpha)-\alpha\, \tilde{B}\gamma_0(\sigma_\alpha)\right]\dd c\\
&=& h\left[ \sum_{i,j=0}^{s-1} \pi_i(\sigma_\alpha)^\top\rho_{ij}(\sigma_\alpha)\gamma_j(\sigma_\alpha)- \alpha\,\pi_0(\sigma_\alpha)^\top \tilde{B}\gamma_0(\sigma_\alpha)\right] ~=~0,
\end{eqnarray*}
provided that (\ref{alfa}) holds true. The statement is completed by observing that the numerator is $O(h^{2s})$, whereas, the denominator is $O(1)$.\QED\bigskip

We now prove that the results of Theorems~\ref{Hcons} and \ref{y1yh} continue to hold for the polynomial (\ref{sigma1alfa}).

\begin{theo}\label{Hconsalfa} For any $\alpha$: $H(\sigma_\alpha(h))=H(\sigma_\alpha(0))$.\end{theo}
\proof Following similar steps as in the proof of Theorem~\ref{Hcons}, ona has:
\begin{eqnarray*}
\lefteqn{H(\sigma_\alpha(h))-H(\sigma_\alpha(0)) = \int_0^h \nabla H(\sigma_\alpha(t))^\top\dot\sigma_\alpha(t)\dd t}\\ 
&=& h\int_0^1 \nabla H(\sigma_\alpha(ch))^\top\dot\sigma_\alpha(ch)\dd c~=~ h\int_0^1 \nabla H(\sigma_\alpha(ch))^\top \sum_{i,j=0}^{s-1} P_i(c)\rho_{ij}(\sigma_\alpha) \gamma_j(\sigma_\alpha)\dd c\\
&&- h\alpha\left[\int_0^1\nabla H(\sigma_\alpha(ch)) \dd c\right]^\top \tilde{B}\gamma_0(\gamma_\alpha)\\
 &=&h \underbrace{\sum_{i,j=0}^{s-1} \gamma_i(\sigma_\alpha)^\top\rho_{ij}(\sigma_\alpha) \gamma_j(\sigma_\alpha)}_{=0} - h\/\alpha\,\gamma_0(\sigma_\alpha)^\top \tilde{B}\gamma_0(\sigma_\alpha)
~=~0,
\end{eqnarray*}
due to the fact that $\tilde{B}$ is skew-symmetric, independently of the considered value of the parameter $\alpha$.\,\QED\bigskip

\begin{theo}\label{y1yhalfa} Assume that the parameter $\alpha$ in (\ref{sigma1alfa}) is chosen according to (\ref{alfa}). Then, $$\sigma_\alpha(h)-y(h)=O(h^{2s+1}).$$\end{theo}\proof Repeating similar steps as those in the proof of Theorem~\ref{y1yh} (and using the same notation), and taking into account (\ref{sigma1alfa}), one arrives at:\footnote{For sake of brevity, we skip here some of the intermediate passages, which are identical to those used in the proof of Theorem~\ref{y1yh}.}
\begin{eqnarray*}
\lefteqn{\sigma_\alpha(h)-y(h)}\\
&=& -h\left[\sum_{i,j\ge0} \Psi_i(\sigma_\alpha)\rho_{ij}(\sigma_\alpha)\gamma_j(\sigma_\alpha) -\sum_{i,j=0}^{s-1} \Psi_i(\sigma_\alpha)\rho_{ij}(\sigma_\alpha)\gamma_j(\sigma_\alpha) +\alpha \Psi_0(\sigma_\alpha)\tilde{B}\gamma_0(\sigma_\alpha)\right]\\
&=&-h\left[ \sum_{i=0}^{s-1}\sum_{j\ge s}  \underbrace{\Psi_i(\sigma_\alpha)\rho_{ij}(\sigma_\alpha)}_{=O(h^j)}\gamma_j(\sigma_\alpha) +  \sum_{i\ge s}\sum_{j=0}^{s-1}  \Psi_i(\sigma_\alpha)\underbrace{\rho_{ij}(\sigma_\alpha)\gamma_j(\sigma_\alpha)}_{=O(h^i)}\right.\\
&&\left. +  \sum_{i,j\ge s}  \Psi_i(\sigma_\alpha)\rho_{ij}(\sigma_\alpha)\gamma_j(\sigma_\alpha)+\underbrace{\alpha \Psi_0(\sigma_\alpha)\tilde{B}\gamma_0(\sigma_\alpha)}_{=O(h^{2s})}
 \right] ~=~O(h^{2s+1}).\,\QED
 \end{eqnarray*}\bigskip

\begin{rem}\label{CasC00}
We observe that the modified polynomial $\sigma_\alpha$ in (\ref{sigma1alfa}) is the solution of the approximate perturbed ODE-IVPs:
\begin{equation}\label{odeivp_salfa}
\dot\sigma_\alpha(ch) = \left[\left(B(\sigma_\alpha(ch))+\alpha \tilde{B}\right)\left[\nabla H(\sigma_\alpha(ch))\right]_s\right]_s, \qquad c\in[0,1], \qquad \sigma_\alpha(0)=y_0,
\end{equation}
where the parameter $\alpha$ is such that $C(\sigma_\alpha(h))=C(\sigma_\alpha(0))$. Clearly, when $\alpha=0$, one recovers the problem (\ref{odeivp_s}) defining $\sigma$.
\end{rem}

We end this section by sketching the case when we have $r$ independent Casimirs, so that $C:\RR^m\rightarrow \RR^r$. In such a case, the notation introduced above formally still holds true, with the following differences:
\begin{itemize}

\item the Fourier coefficients (see (\ref{piy})) $\pi_i(u_\alpha)\in\RR^{m\times r}$, $i=0,\dots,s-1$;

\item the polynomial (\ref{sigma1alfa}) now becomes
\begin{equation}\label{sigma1alfar}
\dot \sigma_\alpha(ch) =   \sum_{i,j=0}^{s-1} P_i(c)\rho_{ij}(\sigma_\alpha) \gamma_j(\sigma_\alpha) - \sum_{\ell=1}^r\alpha_\ell \tilde{B}_\ell\gamma_0(\sigma_\alpha), \qquad c\in[0,1], \qquad \sigma_\alpha(0)=y_0, 
\end{equation}
having set $\alpha=\left(\alpha_1,\,\dots,\,\alpha_r\right)^\top$ and with $\tilde{B}_i^\top=-\tilde{B}_i$, $i=1,\dots,r$, arbitrary skew-symmetric matrices such that  
\begin{equation}\label{M}
M:= \left[ \pi_0(\sigma_\alpha)^\top \tilde{B}_1\gamma_0(\sigma_\alpha),\,\dots,\,\pi_0(\sigma_\alpha)^\top \tilde{B}_r\gamma_0(\sigma_\alpha)\right]\in\RR^{r\times r}
\end{equation}
is nonsingular;
\item the vector $\alpha$, providing the conservation of all Casimirs, is given by (compare with (\ref{alfa}))
\begin{equation}\label{alfar} 
\alpha = M^{-1} \sum_{i,j=0}^{s-1} \pi_i(\sigma_\alpha)^\top \rho_{ij}(\sigma_\alpha)\gamma_j(\sigma_\alpha).
\end{equation}
 
\end{itemize}

\section{Discretization}\label{discr}
The procedure (\ref{sigma1})--(\ref{sigmah}) described in the previous section is not yet a {\em ready to use} numerical method. In fact, in order for this to happen, the integrals $\gamma_j(\sigma), \rho_{ij}(\sigma)$, $i,j=0,\dots,s-1$, defined in (\ref{gammaro}) need to be conveniently computed or approximated. For this purpose, as it has been done in the case of HBVMs \cite{BIT2010}, we shall use a Gauss-Legendre quadrature formula of order $2k$, i.e., the interpolatory quadrature rule based at the zeros of $P_k(c)$, with abscissae and weights $(c_i,b_i)$, for a convenient value $k\ge s$. In so doing, we shall in general obtain a new polynomial approximation $u\in\Pi_s$, in place of $\sigma$ as defined in (\ref{sigma1})--(\ref{sigma}):
\begin{eqnarray}\nonumber
\dot u(ch) =\sum_{i,j=0}^{s-1}P_i(c)\hat\rho_{ij}(u)\hat\gamma_j(u),&&
u(ch) = y_0 + h\sum_{i,j=0}^{s-1}\int_0^c P_i(x)\dd x\hat\rho_{ij}(u)\hat\gamma_j(u),
\qquad c\in[0,1],\\[-3mm] \label{u} \\[-3mm]\nonumber
\hat\gamma_j(u) = \sum_{\ell=1}^k b_\ell P_j(c_\ell)\nabla H(u(c_\ell h)),&&
\hat\rho_{ij}(u) = \sum_{\ell=1}^k b_\ell P_i(c_\ell)P_j(c_\ell)B(u(c_\ell h)),\quad i,j=0,\dots,s-1.
\end{eqnarray}
Consequently, the new approximation to $y(h)$ will be given by 
\begin{equation}\label{uh}
y_1:=u(h) = y_0 + h\sum_{j=0}^{s-1}\hat\rho_{0j}(u) \hat\gamma_j(u)
\equiv y_0 + h\sum_{j=0}^{s-1} \sum_{\ell=1}^k b_\ell P_j(c_\ell)B(u(c_\ell h)) \sum_{\ell=1}^k b_\ell P_j(c_\ell) \nabla H(u(c_\ell h)),
\end{equation}
which is the discrete counterpart of (\ref{sigmah}).  

It is worth mentioning that, in a similar way as done in Section~\ref{intersigma} for the polynomial $\sigma$, for $u$ one obtains, by virtue of (\ref{u}):
\begin{eqnarray*}
\dot u(ch) &=& \sum_{i,j=0}^{s-1}P_i(c)\hat\rho_{ij}(u)\hat\gamma_j(u)\\
&=& \sum_{i=0}^{s-1} P_i(c) \sum_{j=0}^{s-1} \sum_{\ell=1}^k  b_\ell P_i(c_\ell) B(u(c_\ell h))P_j(c_\ell) \sum_{\ell_1=1}^k b_{\ell_1}P_j(c_{\ell_1})\nabla H(u(c_{\ell_1}h))\\
&=& \sum_{i=0}^{s-1} P_i(c)  \sum_{\ell=1}^k b_\ell P_i(c_\ell) B(u(c_\ell h))\sum_{j=0}^{s-1}P_j(c_\ell) \sum_{\ell_1=1}^k b_{\ell_1} P_j(c_{\ell_1})\nabla H(u(c_{\ell_1}h))\\
&=:& \sum_{i=0}^{s-1} P_i(c)  \sum_{\ell=1}^k b_\ell P_i(c_\ell) B(u(c_\ell h))\, [\nabla H(u(c_\ell h))]_s^{(2k)}\\
&\equiv&\left[ B(u(ch))\,[\nabla H(u(ch)]_s^{(2k)}\right]_s^{(2k)},
\end{eqnarray*}
having set $[\cdot]_s^{(2k)}$ the approximate best approximation in {\red $\Pi_{s-1}$} obtained by using a quadrature of order $2k$ for approximating the involved integrals.\footnote{I.e., $\lim_{p\rightarrow\infty} [\cdot]_s^{(p)} = [\cdot]_s$.} Consequently (compare with (\ref{odeivp_s})), the polynomial $u$ is the solution of the initial value problem:
\begin{equation}\label{odeivp_s1}
\dot u(ch) = \left[B(u(ch))\left[\nabla H(u(ch))\right]_s^{(2k)}\right]_s^{(2k)}, \qquad c\in[0,1], \qquad u(0)=y_0.
\end{equation}

\begin{rem}\label{altri}
We observe that the polynomial approximation defined by the problem
$$\dot u(ch) = \left[B(u(ch))\left[\nabla H(u(ch))\right]_s\right]_s^{(2s)}, \qquad c\in[0,1], \qquad u(0)=y_0,$$
corresponds to that provided by the methods in \cite{CH2011}, {\red when the Gauss-Legendre abscissae are used, and to the methods in \cite[Definition\,3.2]{MHW2022}, upon selecting the derivative space as $\Pi_{s-1}$}. Similarly, some of the methods in \cite{BCMR2012} provide the approximation
$$\dot u(ch) = \left[B(u(ch))\left[\nabla H(u(ch))\right]_s^{(2k)}\right]_s^{(2s)}, \qquad c\in[0,1], \qquad u(0)=y_0.$$
\end{rem}

\begin{rem}\label{hbvmks}
When in (\ref{odeivp_s1}) $B(\sigma(ch))\equiv J$, a constant skew-symmetric matrix, we derive the polynomial approximation provided by a HBVM$(k,s)$ method:
\begin{equation}\label{uks}
u(ch) = y_0 + h\int_0^c \left[J\nabla H(u(\tau h))\right]_s^{(2k)}\dd\tau, \qquad c\in[0,1].
\end{equation}
\end{rem}

\subsection{Analysis}\label{analu} As was done in Section~\ref{analsigma} for the continuous procedure, let us now analyze the fully discrete method (\ref{u})--(\ref{uh}). To begin with, the following straightforward result holds true.

\begin{theo}\label{esatto}If 
\begin{equation}\label{klarge}
B\in\Pi_\mu,\quad H\in\Pi_\nu, \qquad \mbox{with}\qquad \mu\le \frac{2k+1}s-2, \quad \nu\le \frac{2k}s.
\end{equation}
Then (see (\ref{u}) and (\ref{gammaro})),
\begin{equation}\label{exact}
\hat\rho_{ij}(u) = \rho_{ij}(u), \qquad \hat\gamma_j(u)=\gamma_j(u), \qquad \forall i,j=0,\dots,s-1,
\end{equation}
and, consequently, with reference to (\ref{u}) and (\ref{sigma}), one has $u\equiv \sigma$.\end{theo}
\proof In fact, if $B$ is a polynomial of degree $\mu$ and $H$ a polynomial of degree $\nu$, the integrand defining $\rho_{ij}(u)$ has at most degree $\mu s +2s-2$, whereas that defining $\gamma_j(u)$ has at most degree $\nu s-1$. Consequently, these degrees do not exceed $2k-1$, when (\ref{klarge}) holds true. As a result, the quadrature is exact, so that (\ref{exact}) is valid and, therefore, $u\equiv\sigma$.\QED\bigskip

Consequently, when (\ref{klarge}) hold true, the method is energy-conserving and has order $2s$, as stated by Theorems~\ref{Hcons} and \ref{y1yh}, respectively.  

Concerning energy conservation, the following additional result holds true, in the case where only $H$ is a polynomial.

\begin{theo}\label{esattoH}If 
\begin{equation}\label{klarge1}
H\in\Pi_\nu, \qquad \mbox{with}\qquad \nu\le \frac{2k}s,
\end{equation}
then $H(y_1)=H(y_0)$.
\end{theo}
\proof In fact, in such a case $\gamma_j(u)=\hat\gamma_j(u)$, ~$j=0,\dots,s-1$, and the proof of Theorem~\ref{Hcons} continues formally to hold, upon replacing $\sigma$ with $u$, and $\rho_{ij}$ with $\hat\rho_{ij}$,
due to the fact that (compare with (\ref{roij}))
\begin{equation}\label{hroij}
\hat\rho_{ij}(u) = \hat\rho_{ji}(u) = -\hat\rho_{ij}(u)^\top, \qquad \forall i,j=0,1,\dots,s-1.\,\QED
\end{equation}

When (\ref{exact}) does not hold true, there is a quadrature error that, upon regularity assumptions, can be easily seen to be given by (see (\ref{gammaro})):
\begin{eqnarray}\label{notexact}
\hat\rho_{ij}(u) - \rho_{ij}(u) &=& \chi_{ij}(h) ~\equiv~ O(h^{2k-i-j}), \\[2mm] \nonumber
\hat\gamma_j(u)-\gamma_j(u) &=& \Delta_j(h)~\equiv~ O(h^{2k-j}), \qquad \forall i,j=0,\dots,s-1.
\end{eqnarray}
Nonetheless, also in this case it is straightforward to verify that  (compare with (\ref{gammaj})),
\begin{equation}\label{hgammaj}
\forall k\ge s: \qquad \hat\gamma_j(u)=O(h^j), \qquad \hat\rho_{ij}(u)=O(h^{|i-j|}).\qquad \forall i,j=0,1,\dots,s-1.
\end{equation}

Consequently, with reference to the approximation $y_1$ defined in (\ref{uh}), the following result is easily obtained,   when (\ref{klarge1}) is not valid.

\begin{theo}\label{Happr} $\forall k\ge s:~H(y_1)=H(y_0)+O(h^{2k+1}).$\end{theo}
\proof In fact, using arguments similar to those used in the proof of Theorem~\ref{Hcons}, one has, by taking into account (\ref{hroij})--(\ref{hgammaj}):
\begin{eqnarray*}
\lefteqn{H(y_1)-H(y_0)~=~ H(u(h))-H(u(0)) = \int_0^h \nabla H(u(t))^\top\dot u(t)\dd t}\\ 
&=& h\int_0^1 \nabla H(u(ch))^\top\dot u(ch)\dd c~=~ h\int_0^1 \nabla H(u(ch))^\top \sum_{i,j=0}^{s-1} P_i(c)\hat\rho_{ij}(u) \hat\gamma_j(u)\dd c\\
 &=&h \sum_{i,j=0}^{s-1} \underbrace{\left[\int_0^1 P_i(c)\nabla H(u(ch))\dd c\right]^\top}_{=\gamma_i(u)^\top}\hat\rho_{ij}(u) \left[\gamma_j(u)+\Delta_j(h)\right]  \\
 &=& h \sum_{i,j=0}^{s-1} \underbrace{\gamma_i(u)^\top\hat\rho_{ij}(u) \gamma_j(u)}_{=0}~+~h \sum_{i,j=0}^{s-1} \gamma_i(u)^\top\hat\rho_{ij}(u) \Delta_j(h)\\
&=&h \left[\sum_{j=0}^{s-1} \sum_{i=j}^{s-1} \underbrace{\gamma_i(u)^\top\hat\rho_{ij}(u)}_{\red =O(h^{2i-j})}\Delta_j(h) +
\sum_{j=0}^{s-1} \sum_{i=0}^{j-1} \gamma_i(u)^\top\underbrace{\hat\rho_{ij}(u)\Delta_j(h)}_{=O(h^{2k-i})}\right] ~=~O(h^{2k+1}).\,\QED
\end{eqnarray*}

\bigskip
Concerning the accuracy of the approximation (\ref{uh}), the following result, stating that the convergence order of Theorem~\ref{y1yh} is retained, holds true.

\begin{theo}\label{y1uh} $\forall k\ge s:~y_1-y(h)=O(h^{2s+1})$.\end{theo}
\proof By using arguments and notations similar to those used in the proof of Theorem~\ref{y1yh}, one has, by taking into account (\ref{notexact}) and that $k\ge s$:
\begin{eqnarray*}
 \lefteqn{
 y_1-y(h) ~=~u(h)-y(h) = y(h,h,u(h))-y(h,0,u(0)) = \int_0^h \frac{\dd}{\dd t} y(h,t,u(t))\dd t}\\
 &=&\int_0^h \left[\left.\frac{\partial}{\partial \xi} y(h,\xi,u(t))\right|_{\xi=t} + \left.\frac{\partial}{\partial \eta} y(h,t,\eta)\right|_{\eta=u(t)}
 \dot u(t)\right] \dd t\\
 &=& \int_0^h\left[ -\Phi(h,t,{\red u(t)})F(u(t))+\Phi(h,t,{\red u(t)})\dot u(t)\right]\dd t \\
 &=& -h\int_0^1 \Phi(h,ch,{\red u(ch)})\left[ F(u(ch))-\dot u(ch)\right]\dd c\\
 &=&-h\int_0^1 \Phi(h,ch,{\red u(ch)})\left[ B(u(ch))\nabla H(u(ch)) - \sum_{i,j=0}^{s-1} P_i(c)\hat\rho_{ij}(u)\hat\gamma_j(u)\right]\dd c\\
 \end{eqnarray*}\begin{eqnarray*}
 &=&\underbrace{-h\int_0^1 \Phi(h,ch,{\red u(ch)})\left[ \sum_{i,j\ge0} P_i(c)\rho_{ij}(u)\gamma_j(u)- \sum_{i,j=0}^{s-1} P_i(c)\rho_{ij}(u)\gamma_j(u)\right]\dd c}_{=O(h^{2s+1}), \mbox{\small ~from the proof of Theorem~\ref{y1yh}}}\\
 &&-h\int_0^1 \Phi(h,ch,{\red u(ch)})\left[ \sum_{i,j=0}^{s-1} P_i(c)\rho_{ij}(u)\gamma_j(u)- \sum_{i,j=0}^{s-1} P_i(c)\hat\rho_{ij}(u)\hat\gamma_j(u)\right]\dd c\\
 &=& O(h^{2s+1}) -h\sum_{i,j=0}^{s-1}\underbrace{\int_0^1 P_i(c) \Phi(h,ch,{\red u(ch)})\dd c}_{=\Psi_i(u)}\left[ \rho_{ij}(u)\gamma_j(u)-\hat\rho_{ij}(u)\hat\gamma_j(u)\right]\\
 &=& O(h^{2s+1}) -h\sum_{i,j=0}^{s-1} \Psi_i(u)\left[ \rho_{ij}(u)\gamma_j(u)-\left(\rho_{ij}(u)+\chi_{ij}(h)\right)\left(\gamma_j(u)+\Delta_j(h)\right)\right]\\
 &=& O(h^{2s+1}) +h\sum_{i,j=0}^{s-1} \Psi_i(u)\left[ \rho_{ij}(u)\Delta_j(h) +\underbrace{\chi_{ij}(h)\gamma_j(u)}_{=O(h^{2k-i})}+\underbrace{\chi_{ij}(h)\Delta_j(h)}_{=O(h^{4k-2j-i})}\right]\\
 &=&O(h^{2s+1}) + O(h^{2k+1}) + h\sum_{i,j=0}^{s-1} \Psi_i(u)\rho_{ij}(u)\Delta_j(h)
 ~=~O(h^{2s+1}) + h\sum_{i,j=0}^{s-1} \Psi_i(u)\rho_{ij}(u)\Delta_j(h).
 \end{eqnarray*}
 Concerning the latter sum, one has:
 \begin{eqnarray*}
 h\sum_{i,j=0}^{s-1} \Psi_i(u)\rho_{ij}(u)\Delta_j(h) 
 &=&h\sum_{i=0}^{s-1}\sum_{j=i}^{s-1} \underbrace{\Psi_i(u)\rho_{ij}(u)}_{=O(h^j)}\Delta_j(h) ~+~h\sum_{i=0}^{s-1}\sum_{j=0}^{i-1} \Psi_i(u)\underbrace{\rho_{ij}(u)\Delta_j(h)}_{=O(h^{2k-2j+i})}\\&=&O(h^{2k+1}) + O(h^{2k+2}) ~=~ O(h^{2k+1}),
 \end{eqnarray*}
 and, consequently, the statement follows.\,\QED
 \medskip
 
 \begin{rem} By taking into account the results of Theorems~\ref{esatto} and \ref{Happr}, it follows that, by choosing $k$ large enough, one obtains either an exact conservation of the Hamiltonian function, in the polynomial case, or a {\em practical} conservation, in the non polynomial case. In fact, as it has been also observed for HBVMs \cite{BIT2012}, in the latter case it is enough to choose $k$ large enough so that the Hamiltonian error falls within the round-off error level of the finite precision arithmetic used in the simulation.
 \end{rem}
 
 It is worth mentioning that a result similar to that of Theorem~\ref{simme} holds true for the fully discrete method.
 
 \begin{theo}\label{simme1} The method (\ref{u})--(\ref{uh}) is symmetric, provided that the abscissae of the quadrature satisfy\,\footnote{This is, in fact, the case, for the Gauss-Legendre quadrature abscissae.}
\begin{equation}
\label{symc}
c_{k-i+1} = 1-c_i, \qquad i=1,\dots,k.
\end{equation}
\end{theo} 
\proof  
Symmetry of a given one-step method $y_1=\Phi_h(y_0)$ applied to an initial value problem $y'=f(y)$ with $y(0)=y_0$, means that $\Phi_h^{-1}=\Phi_{-h}$ that is, applying the method to the state vector $y_1$, but with the direction of time reversed, yields the initial state vector $y_0$, independently of the choice of the initial value $y_0$. In our context, with reference to (\ref{u})--(\ref{uh}), $\Phi_h$ is defined by
\begin{equation}
\label{Phi_h}
y_1 = y_0 + h\sum_{j=0}^{s-1}\hat\rho_{0j} \hat\gamma_j,
\end{equation}
where $\hat\rho_{ij}$ and $\hat\gamma_j$ are the solutions of the nonlinear system
\begin{eqnarray}\nonumber
\hat\gamma_j &=& \sum_{\ell=1}^k b_\ell P_j(c_\ell)\nabla H\hspace{-1mm}\left( y_0 + h\sum_{\mu,\nu=0}^{s-1}\int_0^{c_\ell} P_\mu(x)\dd x\hat\rho_{\mu\nu}\hat\gamma_\nu \right), \\  \nonumber
\hat\rho_{ij} &=& \sum_{\ell=1}^k b_\ell P_i(c_\ell)P_j(c_\ell)B\hspace{-1mm}\left( y_0 + h\sum_{\mu,\nu=0}^{s-1}\int_0^{c_\ell} P_\mu(x)\dd x\hat\rho_{\mu\nu}\hat\gamma_\nu \right),\\  \label{nonlinPhi_h}
&&i,j=0,\dots,s-1.
\end{eqnarray}
We can obtain the explicit formulation of  $\bar y_0=\Phi_{-h}(y_1)$ by  introducing in (\ref{Phi_h})--(\ref{nonlinPhi_h}) the following substitutions:
$y_1$ in place of $y_0$ and  $-h$ in place of $h$. In so doing, we arrive at the method defined as
\begin{equation}
\label{bu}
\bar y_0 = y_1 - h\sum_{j=0}^{s-1}\bar\rho_{0j} \bar\gamma_j,
\end{equation}
where the unknown quantities $\bar\rho_{ij}$ and $\bar\gamma_j$ satisfy the following nonlinear system 
\begin{eqnarray}\nonumber
\bar\gamma_j &=& \sum_{\ell=1}^k b_\ell P_j(c_\ell)\nabla H\hspace{-1mm}\left(y_1 - h\sum_{\mu,\nu=0}^{s-1}\int_0^{c_\ell} P_\mu(x)\dd x\bar\rho_{\mu\nu}\bar\gamma_\nu\right), \\ \nonumber
\displaystyle \bar\rho_{ij} &=&  \sum_{\ell=1}^k b_\ell P_i(c_\ell)P_j(c_\ell)B\hspace{-1mm}\left(y_1 - h\sum_{\mu,\nu=0}^{s-1}\int_0^{c_\ell} P_\mu(x)\dd x\bar\rho_{\mu\nu}\bar\gamma_\nu\right), \\
&&i,j=0,\dots,s-1,  \label{nonlinPhi_{-h}}
\end{eqnarray}
and we want to show that $\bar y_0=y_0$. To this end, we introduce the following variables:
$$
\gamma_j^\ast := (-1)^j \hat \gamma_j, \quad \rho_{ij}^\ast := (-1)^{i+j} \hat \rho_{ij}, \quad i,j=0,\dots,s-1.
$$
Exploiting the symmetry property of Legendre polynomials, $(-1)^jP_j(c)=P_j(1-c)$, from the first equation in (\ref{nonlinPhi_h}) we get
\begin{eqnarray*}
\gamma_j^\ast &=& \displaystyle \sum_{\ell=1}^k b_\ell (-1)^jP_j(c_\ell)\nabla H\hspace{-1mm}\left( y_0 + h\sum_{\mu,\nu=0}^{s-1}\left(\int_0^1 P_\mu(x){\red \dd x} - \int_{c_\ell}^1 P_\mu(x)\dd x \right)\hat\rho_{\mu\nu}\hat\gamma_\nu \right) \\
&=& \sum_{\ell=1}^k b_\ell P_j(1-c_\ell)\nabla H\hspace{-1mm}\left( y_1 - h\sum_{\mu,\nu=0}^{s-1} \int_{c_\ell}^1 P_\mu(x)\dd x \hat\rho_{\mu\nu}\hat\gamma_\nu \right). 
\end{eqnarray*}
Introducing the change of variables $\tau=1-x$, transforms the latter integral as
$$
\int_{c_\ell}^1 P_\mu(x)\dd x = -\int_{1-c_\ell}^0 P_\mu(1-\tau)\dd \tau = (-1)^\mu \int_0^{1-c_\ell} P_\mu(\tau)\dd \tau. 
$$
Exploiting the symmetry assumption (\ref{symc}), which in turn implies symmetric weights, $b_{k-i+1} = b_i$, $i=1,\dots,k$, we finally get
\begin{eqnarray*}
\gamma_j^\ast &=&  \displaystyle \sum_{\ell=1}^k b_{k-\ell+1} P_j(c_{k-\ell+1})\nabla H\hspace{-1mm}\left( y_1 - h\sum_{\mu,\nu=0}^{s-1} \int_0^{c_{k-\ell+1}} P_\mu(x)\dd x  (-1)^\mu (-1)^\nu \hat\rho_{\mu\nu} (-1)^\nu \hat \gamma_\nu \right)\\ 
&=&  \sum_{\ell=1}^k b_{\ell} P_j(c_{\ell})\nabla H\hspace{-1mm}\left( y_1 - h\sum_{\mu,\nu=0}^{s-1} \int_0^{c_{\ell}} P_\mu(x)\dd x   \rho^\ast_{\mu\nu}  \gamma^\ast_\nu \right).
\end{eqnarray*}
The same flow of computation may be employed on the second equation in (\ref{nonlinPhi_h}) to see that
$$
\displaystyle \rho^\ast_{ij} =  \sum_{\ell=1}^k b_\ell P_i(c_\ell)P_j(c_\ell)B\hspace{-1mm}\left(y_1 - h\sum_{\mu,\nu=0}^{s-1}\int_0^{c_\ell} P_\mu(x)\dd x\rho^\ast_{\mu\nu}\gamma^\ast_\nu\right).
$$ 
We then realize that $\gamma_j^\ast$ and   $\rho^\ast_{ij}$ satisfy the very same nonlinear system (\ref{nonlinPhi_{-h}}) governing the quantities $\bar \gamma_j$ and   $\bar \rho_{ij}$. Thus we may conclude that
$$
\bar \gamma_j= (-1)^j \hat \gamma_j, \quad \bar \rho_{ij}= (-1)^{i+j} \hat \rho_{ij}, \quad i,j=0,\dots,s-1,
$$
and hence from (\ref{bu}),
$$
\bar y_0 = y_1 - h\sum_{j=0}^{s-1}(-1)^j \hat \rho_{0j} (-1)^j \hat \gamma_j =  y_1 - h\sum_{j=0}^{s-1} \hat \rho_{0j}  \hat \gamma_j =y_0.
$$\,\QED\bigskip

In the limit case when $k\rightarrow \infty$, since the Gauss-Legendre quadrature formule are convergent, we are led back to the symmetry property of the original non-discretized procedure (\ref{sigma1})--(\ref{sigmah}) that we anticipated in Theorem \ref{simme}.

We conclude this section by emphasizing that, when problem (\ref{poisson}) is in the form (\ref{Ham}), then (see (\ref{u})) $\hat\rho_{ij}(\sigma) = {\red \delta_{ij}}J$ and, consequently, the polynomial   approximation $u$ becomes
 $$u(ch) = y_0+h\sum_{j=0}^{s-1} \int_0^c P_j(x)\dd x \sum_{{\red i=1}}^k b_i P_j(c_i)J\nabla H(u(c_i h)), \qquad c\in[0,1],$$
 which is equivalent to (\ref{uks}). As anticipated in Remark~\ref{hbvmks}, this equation (see, e.g., the monograph \cite{LIMbook2016} or the review paper \cite[Section\,2.2]{BI2018}) defines a {\em Hamiltonian Boundary Value Method with parameters $k$ and $s$, in short HBVM$(k,s)$}. Moreover,  also Theorems~\ref{esatto}--\ref{y1uh} exactly describe, in the case of problem (\ref{Ham}), their conservation and accuracy properties. For this reason, the methods here presented can be regarded as a generalization of HBVMs for solving Poisson problems and, therefore,  this fact, motivates the following definition.
 
\begin{defi}\label{phbvm}  We shall refer to the method (\ref{u})--(\ref{uh}) as to a {\em PHBVM$(k,s)$ method}. \end{defi}

\subsection{Conservation of Casimirs}\label{cassec_d}
This section is devoted to the conservation of Casimirs for the discrete PHBVM$(k,s)$ method (\ref{u})--(\ref{uh}), following similar steps as those in Section~\ref{cassec}. For this purpose, it is convenient to define the approximate Fourier coefficients, for $k\ge s$, which we assume hereafter:
\begin{equation}\label{pisig}
\hat\pi_i(\sigma) = \sum_{\ell=1}^k b_\ell P_i(c_\ell)\nabla C(\sigma(c_\ell h)) = O(h^j), \qquad i=0,\dots,s-1.
\end{equation}
The following straightforward result is reported without proof.

\begin{lem}\label{pik} Assume that $\sigma\in\Pi_s$. Then, with reference to (\ref{piy}), for all $i=0,\dots,s-1$, one has:
\begin{eqnarray*}
\hat\pi_i(\sigma) &=& \pi_i(\sigma), \qquad \mbox{if} \qquad C\in\Pi_\nu,\quad \nu\le 2k/s,\\[2mm]
\hat\pi_i(\sigma) &=& \pi_i(\sigma) - \theta_i(h), \qquad \theta_i(h)=O(h^{2k-i}), \qquad \mbox{otherwise.}
\end{eqnarray*}
\end{lem}

Next, let us consider the following perturbed polynomial, in place of the polynomial $u$ in (\ref{u}),
\begin{equation}\label{u1alfa}
\dot u_\alpha(ch) =\sum_{i,j=0}^{s-1}P_i(c)\hat\rho_{ij}(u_\alpha)\hat\gamma_j(u_\alpha) -\alpha \tilde{B}\hat\gamma_0(u_\alpha), \qquad c\in[0,1], \qquad u_\alpha(0)=y_0,
\end{equation}
with $\tilde{B}^\top=-\tilde{B}\ne O$ an arbitrary matrix, and (compare with (\ref{alfa}))
\begin{equation}\label{alfad} 
\alpha = \frac{\sum_{i,j=0}^{s-1} \hat\pi_i(u_\alpha)^\top \hat\rho_{ij}(u_\alpha)\hat\gamma_j(u_\alpha)}{\hat\pi_0(u_\alpha)^\top \tilde{B}\hat\gamma_0(u_\alpha)} = O(h^{2s}).
\end{equation}
 
\begin{theo}\label{esatto1}If 
\begin{equation}\label{klarged}
B\in\Pi_\mu,\quad C,H\in\Pi_\nu, \qquad \mbox{with}\qquad \mu\le \frac{2k+1}s-2, \quad \nu\le \frac{2k}s.
\end{equation}
Then (see (\ref{u}), (\ref{gammaro})), (\ref{piy}), and (\ref{pisig})
\begin{equation}\label{exactd}
\hat\rho_{ij}(u_\alpha) = \rho_{ij}(u_\alpha), \qquad \hat\gamma_j(u_\alpha)=\gamma_j(u_\alpha), \qquad \hat\pi_i(u_\alpha)=\pi_i(u_\alpha),\qquad \forall i,j=0,\dots,s-1,
\end{equation}
and, consequently, with reference to (\ref{u1alfa}) and (\ref{sigma1alfa}), one has $u_\alpha\equiv \sigma_\alpha$.\end{theo}

Differently, by setting $y_1:=u_\alpha(h)$ the new approximation, Theorems~\ref{esattoH}, \ref{Happr}, \ref{y1uh}, and \ref{simme1} continue formally to hold. Moreover, the following result holds true.

\begin{theo}\label{Cappr} With reference to (\ref{u1alfa})--(\ref{alfad}), one has: $$
C(y_1)-C(y_0)=\left\{\begin{array}{cc}
0, &\mbox{if}\quad C\in\Pi_\nu, \quad\mbox{with}\quad  \nu\le 2k/s,\\[2mm]
O(h^{2k+1}), &\mbox{otherwise}.\end{array}\right.$$\end{theo}
\proof By taking into account Lemma~\ref{pik}, one obtains:
\begin{eqnarray*}
\lefteqn{C(y_1)-C(y_0) ~=~ h\int_0^1 \nabla C(u_\alpha(ch))^\top \dot u_\alpha(ch)\dd c}\\
&=& h\int_0^1 \nabla C(u_\alpha(ch))^\top \left[ \sum_{i,j=0}^{s-1} P_i(c)\hat\rho_{ij}(u_\alpha)\hat\gamma_j(u_\alpha) - \alpha\,\tilde{B}\hat\gamma_0(u_\alpha)\right]\dd c\\
&=& h\left[\sum_{i,j=0}^{s-1} \left(\int_0^1 \nabla C(u_\alpha(ch))P_i(c)\dd c\right)^\top \hat\rho_{ij}(u_\alpha)\hat\gamma_j(u_\alpha) - \alpha\,\left(\int_0^1 \nabla C(u_\alpha(ch))P_i(c)\dd c\right)^\top \tilde{B}\hat\gamma_0(u_\alpha)\right]\dd c\\
&=& h\left[\sum_{i,j=0}^{s-1} \pi_i(u_\alpha)^\top \hat\rho_{ij}(u_\alpha)\hat\gamma_j(u_\alpha) - \alpha\,\pi_0(u_\alpha)^\top \tilde{B}\hat\gamma_0(u_\alpha)\right]~=:~(*).
\end{eqnarray*}
In case $C\in\Pi_\nu$ with $\nu\le 2k/s$, then, by virtue of Lemma~\ref{pik}, $\pi_i(u_\alpha)=\hat\pi_i(u_\alpha)$ and, consequently, $(*)=0$ because of (\ref{alfad}). Conversely, again from Lemma~\ref{pik}, one has:
\begin{eqnarray*}
(*)&=& h\left[\sum_{i,j=0}^{s-1} \left[\hat\pi_i(u_\alpha)+\theta_i(h)\right]^\top \hat\rho_{ij}(u_\alpha)\hat\gamma_j(u_\alpha) - \alpha\,\left[\hat\pi_0(u_\alpha)+\theta_0(h)\right]^\top \tilde{B}\hat\gamma_0(u_\alpha)\right]\\
&=& h\left[\underbrace{\sum_{i,j=0}^{s-1} \theta_i(h)^\top \hat\rho_{ij}(u_\alpha)\hat\gamma_j(u_\alpha)}_{=O(h^{2k})} - \underbrace{\alpha\,\theta_0(h)^\top \tilde{B}\hat\gamma_0(u_\alpha)}_{=O(h^{2k+2s})}\right]
~=~O(h^{2k+1}).\,\QED
\end{eqnarray*}

\begin{rem} From the arguments here exposed, one deduces that when using a finite precision arithmetic, a (at least) practical conservation of  both the Hamiltonian and the Casimirs can be obtained, by choosing $k$ large enough. \end{rem}

\begin{defi}\label{ephbvm} Following \cite{BS2014}, we name {\em enhanced PHBVM$(k,s)$}, in short\, {\em EPHBVM$(k,s)$}, the method defined by (\ref{u1alfa})--(\ref{alfad}).\end{defi}

Following  Remark~\ref{CasC00}, the polynomial $u_\alpha$ defined by an EPHBVM$(k,s)$ method is the solution of the ODE-IVP (compare with (\ref{odeivp_salfa})):
$$
\dot u_\alpha(ch) = \left[\left(B(u_\alpha(ch))+\alpha \tilde{B}\right)\left[\nabla H(u_\alpha(ch))\right]_s^{(2k)}\right]_s^{(2k)}, \qquad c\in[0,1], \qquad u_\alpha(0)=y_0.
$$
Clearly, when $\alpha=0$ one retrieves the problem  (\ref{odeivp_s1}) defining the polynomial approximation of the PHBVM$(k,s)$ method.

\begin{rem}\label{kggs}
From the results of Theorems~\ref{esattoH}, \ref{Happr}, and \ref{Cappr}, one deduces the clear advantage of choosing values of $k$ suitably larger than $s$, in order to obtain a suitable conservation of non-quadratic Hamiltonians and/or Casimirs. This fact will be duly confirmed in the numerical tests. Moreover, this is not a serious drawback, from the point of view of the computational cost, since, as we shall see in the next section, the discrete problem to be solved will always have (block) dimension $s$, independently of $k$.\footnote{This latter feature is inherited, in turn, from the original methods, HBVM$(k,s)$ and EHBVM$(k,s)$ for Hamiltonian problems.}
\end{rem}

\section{The discrete problem}\label{discrete}
In this section we deal with the efficient solution of the discrete problem generated by the PHBVM$(k,s)$ method (\ref{u}). 
For this purpose, we observe that only the values \,$Y_\ell:=u(c_\ell h)$,\, $\ell=1,\dots,k$, are actually needed. Consequently, (\ref{u}) can be re-written as:
\begin{eqnarray}\label{uell}
Y_\ell &=& y_0 + h\sum_{i,j=0}^{s-1}\int_0^{c_\ell} P_i(x)\dd x\hat\rho_{ij}\hat\gamma_j,
\qquad \ell = 1,\dots,k,\\  \nonumber
\hat\gamma_j &=& \sum_{\ell=1}^k b_\ell P_j(c_\ell)\nabla H(Y_\ell),\qquad
\hat\rho_{ij} ~=~ \sum_{\ell=1}^k b_\ell P_i(c_\ell)P_j(c_\ell)B(Y_\ell),\qquad i,j=0,\dots,s-1,
\end{eqnarray}
where, for sake of brevity, we have omitted the argument $u$ for $\hat\gamma_j$ and $\hat\rho_{ij}$, as was already done in (\ref{Phi_h})--(\ref{nonlinPhi_h}). The equations (\ref{uell}) can be cast in matrix form by defining the block vectors and matrices
\begin{equation}\label{Gg}
Y = \pmatrix{c} Y_1\\ \vdots \\ Y_k\endpmatrix\in\RR^{k\cdot m}, \quad \bfgamma = \pmatrix{c}\hat\gamma_0\\ \vdots\\ \hat\gamma_{s-1}\endpmatrix\in\RR^{s\cdot m}, \quad 
\Gamma = \pmatrix{ccc} \hat\rho_{00} & \dots &\hat\rho_{0,s-1}\\ \vdots & &\vdots\\ \hat\rho_{s-1,0} & \dots &\hat\rho_{s-1,s-1}\endpmatrix\in\RR^{s\cdot m \times s\cdot m},
\end{equation}
and  
\begin{equation}\label{PIO}
\bfe = \pmatrix{c} 1\\ \vdots\\1\endpmatrix \in \RR^k, \quad \Omega = \pmatrix{ccc} b_1\\ &\ddots\\ &&b_k\endpmatrix, \quad \I_s=\left( \int_0^{c_i} P_{j-1}(x)\dd x\right), ~ \P_s=\left( P_{j-1}(c_i)\right) \in\RR^{k\times s}.
\end{equation}
In fact, by also setting hereafter $I_r$ the identity of dimension $r$, we can rewrite (\ref{uell}) as:
\begin{equation}\label{YgG}
Y = \bfe\otimes y_0 + h(\I_s\otimes I_m) \Gamma\bfgamma, \quad
\bfgamma = (\P_s^\top\Omega\otimes I_m)\nabla H(Y), \quad
 \Gamma = (\P_s^\top\Omega\otimes I_m) \B(Y) (\P_s\otimes I_m),
\end{equation}
where
$$\nabla H(Y) = \pmatrix{c} \nabla H(Y_1)\\ \vdots\\ \nabla H(Y_k)\endpmatrix\qquad\mbox{and}\qquad \B(Y) = \pmatrix{ccc} B(Y_1)\\ &\ddots\\ && B(Y_k)\endpmatrix.$$
As is clear, the discrete problem (\ref{YgG}) can be further reformulated in terms of the product of the Fourier coefficients $\hat\rho_{ij}$ and $\hat\gamma_j$. In fact, by setting
\begin{equation}\label{bffi}
\bfphi \equiv \pmatrix{c} \phi_0\\ \vdots \\ \phi_{s-1}\endpmatrix :=\Gamma\bfgamma \qquad \Rightarrow\qquad \phi_i=\sum_{j=0}^{s-1}\hat\rho_{ij}\hat\gamma_j, \quad i=0,\dots,s-1,
\end{equation}
the first equation in (\ref{YgG}) becomes 
\begin{equation}\label{uno}
Y = \bfe\otimes y_0 + h\I_s\otimes I_m\bfphi,
\end{equation}
whereas, multiplying side by side the third by the second gives 
\begin{equation}\label{due}
\bfphi = (\P_s^\top\Omega\otimes I_m)\, \B(Y)\, (\P_s\P_s^\top\Omega\otimes I_m)\nabla H(Y).
\end{equation}
Consequently, by substituting the right-hand side of (\ref{uno}) in the right-hand side of (\ref{due}), provides the new discrete problem:
\begin{eqnarray}\nonumber
\F(\bfphi)&:=&\bfphi - \P_s^\top\Omega\otimes I_m\, \B\left(\bfe\otimes y_0 + h\I_s\otimes I_m \bfphi\right)\, (\P_s\P_s^\top\Omega\otimes I_m)\nabla H\left(\bfe\otimes y_0 + h\I_s\otimes I_m \bfphi\right)\\ &=& \bfzero.\label{fi}
\end{eqnarray}
Moreover, computing the vector $\bfphi$ in (\ref{bffi}) allows us to obtain the new approximation (\ref{uh}) as:
$$y_1 = y_0 + h\phi_0.$$

\medskip
\begin{rem}\label{HBVM3} In case where the problem (\ref{poisson}) is in the form (\ref{Ham}), one has that $\B(Y)=I_s\otimes J$. Consequently,  considering that
\begin{equation}\label{POP}
\P_s^\top\Omega\P_s=I_s,
\end{equation}
the discrete problem (\ref{fi}) reduces to:
$$
\bfphi - \P_s^\top\Omega\otimes J \nabla H\left(\bfe\otimes y_0 + h\I_s\otimes J \bfphi\right) = \bfzero.
$$
This latter problem is exactly that generated by a HBVM$(k,s)$ method applied for solving (\ref{Ham}) \cite{BIT2011}. We observe, however, that while the original HBVM$(k,s)$ method is actually a $k$-stage Runge-Kutta method with Butcher tableau (see (\ref{PIO}))
$$\begin{array}{c|c} \bfc & \I_s\P_s^\top\Omega \\[2mm] \hline \\[-2mm] & \bfb^\top\end{array}, \qquad \bfb = \pmatrix{c} b_1\\ \vdots \\ b_k\endpmatrix, 
 \qquad \bfc = \pmatrix{c} c_1\\ \vdots \\ c_k\endpmatrix,$$ this is no more the case for the gneralization defined by (\ref{YgG}). 
\end{rem}

\begin{rem}\label{Gauss}
In the case $k=s$, one has that $\P_s\P_s^\top\Omega=I_s$. Consequently, (\ref{fi}) becomes, by using the notation (\ref{poissonxi}),
$$\bfphi = (\P_s^\top\Omega\otimes I_m)\, F(\bfe\otimes y_0 + h\I_s\otimes I_m \bfphi).$$
This, in turn, is equivalent  to the application of the $s$-stage Gauss method to the problem (\ref{poisson}) (see, e.g., \cite{BIT2011}).

Instead, in the case $k>s$, the discrete problem (\ref{fi}) is equivalent to the application of the HBVM$(k,s)$ method to the problem (see (\ref{odeivp_s1})),
$$\dot y = B(y)[\nabla H(y)]_s^{(2k)}, \qquad t>0, \qquad y(0)=y_0,$$
in place of (\ref{poisson}). This application , in turn, provides the polynomial approximation (\ref{odeivp_s1}).
\end{rem}

\medskip
We observe that the formulation (\ref{fi}) naturally induces a straightforward iterative procedure for solving the discrete problem,
\begin{eqnarray}\nonumber
\lefteqn{\bfphi^{r+1} =}\\ \nonumber
&& \P_s^\top\Omega\otimes I_m\, \B\left(\bfe\otimes y_0 + h\I_s\otimes I_m \bfphi^r\right)\, (\P_s\P_s^\top\Omega\otimes I_m)\nabla
H\left(\bfe\otimes y_0 + h\I_s\otimes I_m \bfphi^r\right),\\  \label{ite}
&& r=0,1,\dots,
\end{eqnarray}
for which the initial approximation $\bfphi^0=\bfzero$ can be conveniently used. It is also possible to use the simplified Newton iteration for solving (\ref{fi}) which, taking into account (\ref{PIO}), (\ref{POP}), and that (see, e.g., \cite{LIMbook2016})
\begin{equation}\label{Xs}
\P_s^\top\Omega\I_s = X_s := \pmatrix{cccc}
\xi_0 &-\xi_1\\
\xi_1 &0 &\ddots\\
        &\ddots &\ddots &-\xi_{s-1}\\
        &           &\xi_{s-1} &0\endpmatrix, \qquad \xi_i = \left(2\sqrt{|4i^2-1|}\right)^{-1},\quad i=0,\dots,s-1,~~
\end{equation}
takes the form:
\begin{equation}\label{Newt}
        \mbox{solve:}~\left[ I_s\otimes I_m -hX_s\otimes F'(y_0)\right] \bfdelta^r = -\F(\bfphi^r), \qquad \bfphi^{r+1}:=\bfphi^r+\bfdelta^r, \qquad r=0,1,\dots,
\end{equation}
with $F'$ the Jacobian of $F$ (see (\ref{poissonxi})). 
 Nevertheless, the iteration (\ref{Newt}) requires the factorization of a matrix having size $s$ times larger than that of problem (\ref{poisson}), which can be costly, when $s$ and/or $m$ are large. Consequently, it is much more effective to resort to a {\em blended} iteration  for solving (\ref{fi}) (see, e.g., \cite{BIT2011}, we also refer to \cite{BM2002} for a more detailed analysis of {\em blended methods}). In the present case, this latter iteration, considering the matrix $X_s$ defined in (\ref{Xs}), denoting $\sigma(X_s)$ its spectrum, and setting
 \begin{equation}\label{Sigma}
\Lambda := I_m-h\lambda_s F'(y_0)\in\RR^{m\times m}, \qquad\mbox{with}\qquad \lambda_s = \min_{\lambda\in\sigma(X_s)} |\lambda|,
\end{equation}
assumes the form:        
\begin{eqnarray}\label{blend}
\bfeta^r&:=&-\F(\bfphi^r),\qquad\quad \bfeta_1^r~:=~ \lambda_sX_s^{-1}\otimes I_m\, \bfeta^r,\\[2mm] \nonumber
\bfphi^{r+1}&:=&\bfphi^r+I_s\otimes \Lambda^{-1}\left( \bfeta_1^r +I_s\otimes \Lambda^{-1}\left(\bfeta^r-\bfeta_1^r\right)\right), \qquad r=0,1,\dots.
\end{eqnarray}
Consequently, only the matrix $\Lambda$ in (\ref{Sigma}), having the same size as that of the continuous problem (\ref{poisson}), needs to be factored. This fact is a common feature, in the many instances where the blended iteration can be used. For this reason, in such cases, it turns out to be extremely efficient (see, e.g., \cite{BM2004,BMM2006,BM2009,BFCI2014,WMF2017}). 

\medskip
{\blue
\begin{rem}
In the practical use of the methods, it is customary to choose the parameter $k$, related to the order of the quadrature, so that the discretization error falls within the round-off error level. Nevertheless, round-off errors are unavoidable, as are iteration errors in (\ref{ite}) or (\ref{blend}). This may cause a small numerical drift in the invariants,\footnote{\blue Usually, this is a very small drift, of the order of the machine epsilon per step.} even in the  case where the quadrature is exact. This phenomenon has been duly studied in \cite[Chapter~4.3]{LIMbook2016}, where a simple correction procedure is given to avoid this problem. The same procedure can be conveniently used in this setting, too. The reader is referred to the above reference for full details. 
\end{rem}
}

\subsection{Conservation of Casimirs}\label{ephbvm_sec}
In this section we sketch the implementation of EPHBVM$(k,s)$ methods described in Section~\ref{cassec_d}. For this purpose, besides the vector $\bfphi$ defined in (\ref{bffi}), we need to define the block vector
\begin{equation}\label{bfpi}
\hat\bfpi = \pmatrix{c} \hat\pi_0 \\ \vdots \\ \hat\pi_{s-1}\endpmatrix 
\end{equation}
with the approximate Fourier coefficients (\ref{pisig}) of the gradient of the Casimir.\footnote{As before, for sake of brevity we now omit the argument $u_\alpha$ of the approximate Fourier coefficients.} 
In so doing, the discrete problem generated by an EHBVM$(k,s)$ method becomes: 
\begin{eqnarray}\label{fialfa}
\F(\bfphi,\alpha)&:=&\pmatrix{c}
\bfphi - \P_s^\top\Omega\otimes I_m\, \B\left(Y\right)\, (\P_s\P_s^\top\Omega\otimes I_m)\nabla H\left(Y\right)\\[2mm]
\alpha - \frac{\hat\bfpi^\top\bfphi}{\hat\pi_0^\top \tilde{B}\hat\gamma_0 }\endpmatrix~=~ \bfzero,\\ \nonumber
\mbox{with}&&\\ \nonumber
Y&=& \bfe\otimes y_0 + h\I_s\otimes I_m \bfphi - \alpha h\bfc\otimes (\tilde{B}\hat\gamma_0),\\ \nonumber
\hat\gamma_0 &=&\bfb^\top\otimes I_m\,\nabla H(Y),\\ \nonumber
\hat\bfpi &=& \P_s^\top\Omega\otimes I_m\, \nabla C(Y),
\end{eqnarray}
and the new approximation given by
$$y_1 = y_0 + h\left(\phi_0-\alpha\, \tilde{B}\hat\gamma_0\right).$$

We conclude this section by mentioning that, in case of multiple Casimirs, the discrete problem (\ref{fialfa}) can be readily generalized, by considering discrete counterparts of (\ref{M})--(\ref{alfar}).

\section{Numerical tests}\label{num}

In this section we present a couple of numerical tests concerning the solution of Lotka-Volterra problems, with the last one possessing a Casimir. The numerical tests have been carried out on a 3GHz Intel Xeon W10 core computer with 64GB of memory running Matlab 2020a.

\begin{figure}[t]
\centerline{\includegraphics[width=9.5cm]{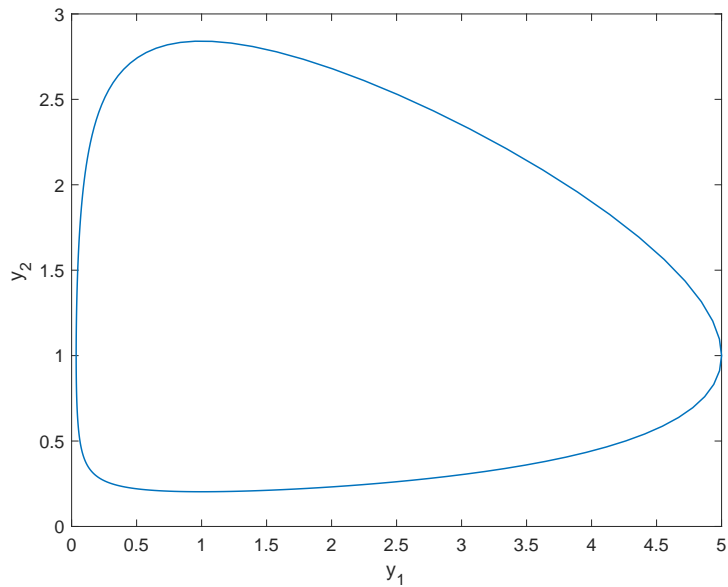}}
\caption{solution of problem (\ref{ex1}).}
\label{fig1}
\end{figure}

\begin{figure}[t]
\centerline{\includegraphics[width=7.25cm]{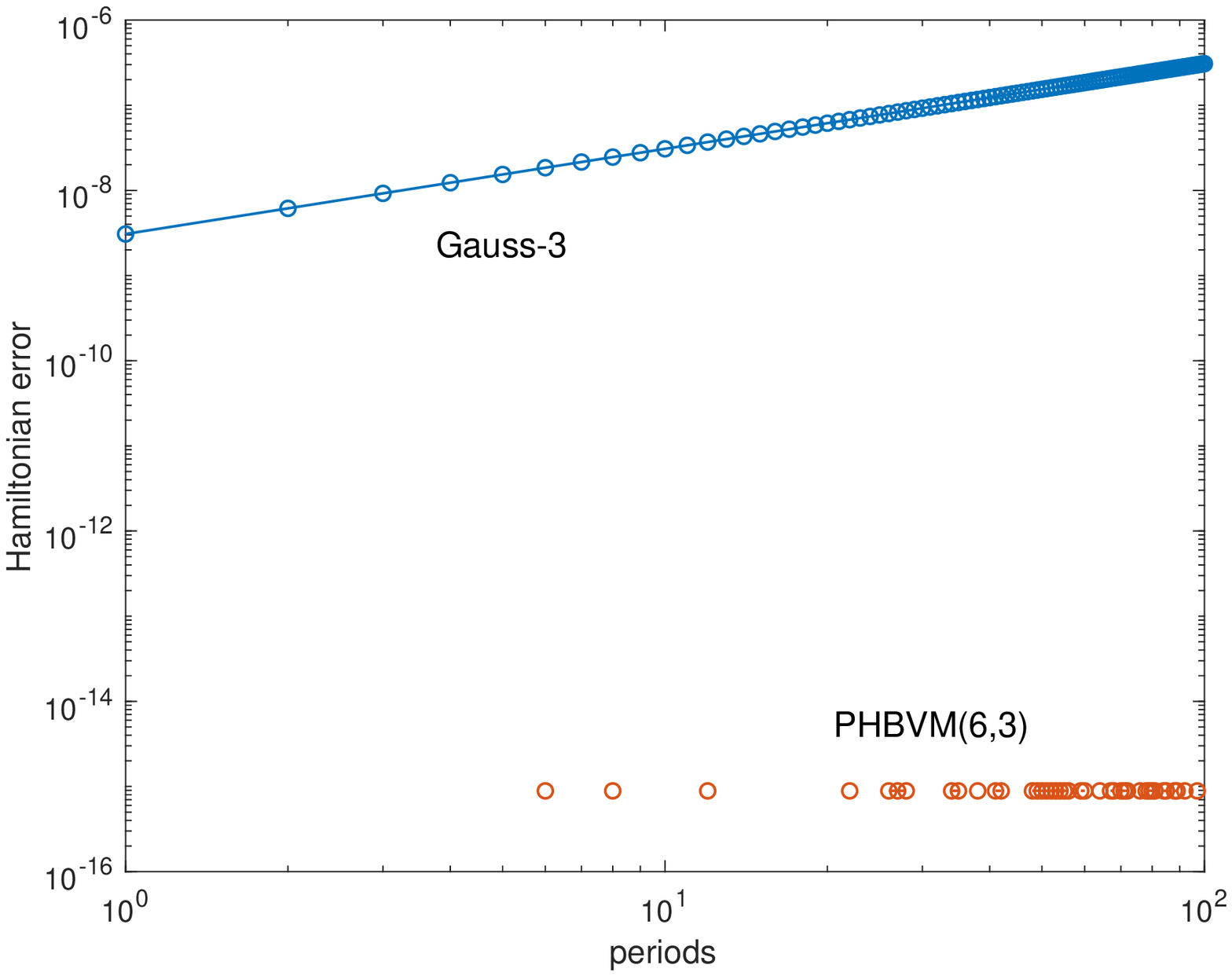}\qquad\includegraphics[width=7.25cm]{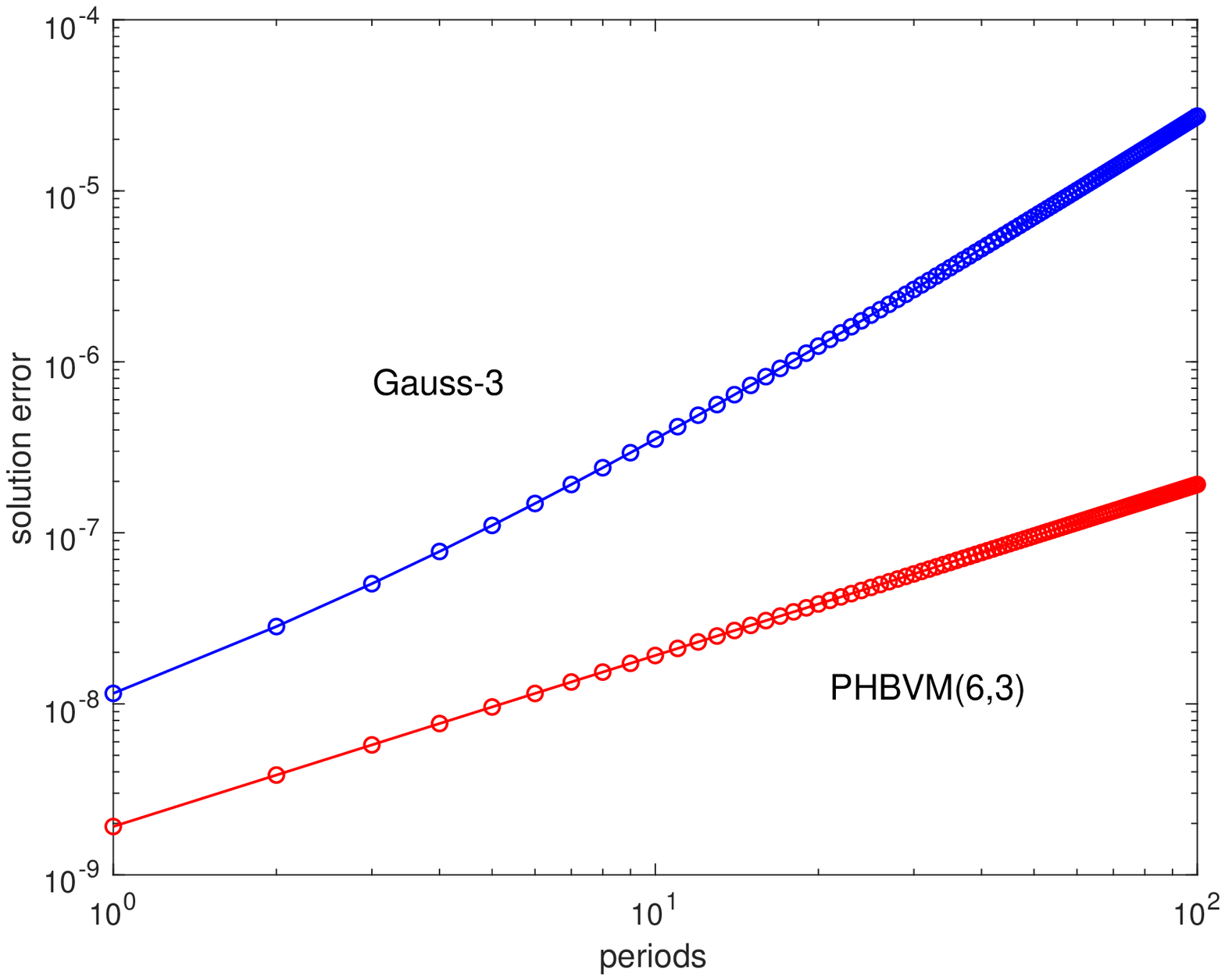}}
\caption{Hamiltonian error (left plot) and solution error (right plot) when solving problem (\ref{ex1}) with time-step $h=T/100$ over 100 periods.}
\label{ex1_err}
\end{figure}

\bigskip
{\bf Example~1.~} We consider the following Lotka-Volterra problem: 
\begin{eqnarray}\nonumber
\dot y &=& \pmatrix{cc} 0 & y_1y_2 \\ -y_1y_2 &0 \endpmatrix \nabla H(y),\\[2mm] \label{ex1}
H(y)   &=& a\left( \ln y_1 -\frac{y_1}{y_1^*}\right) + b\left( \ln y_2 -\frac{y_2}{y_2^*}\right),
\end{eqnarray}
with   
$$a=1, \qquad b=3, \qquad y_1^*=y_2^*=1, \qquad y(0) = \left(5,\,1\right)^\top,$$
whose solution, which is periodic of period $T\approx 4.633434168477889$, is depicted in Figure~\ref{fig1}. At first, we solve the problem on one period with timestep $h=T/n$, by using the following methods:

\begin{itemize}
\item the $s$-stage Gauss method, $s=1,2,3$;
\item the PHBVM$(4,s)$, $s=1,2$, and PHBVM$(6,3)$ methods, which become soon energy-conserving, as the value of $n$ is increased.
\end{itemize}
The obtained results are summarized in Table~\ref{ex1tab}, where we have denoted by $e_y$ and $e_H$ the error in the solution and in the Hamiltonian after one period, respectively. Their numerical rate of convergence is also reported, along with the mean blended iterations (\ref{Sigma})--(\ref{blend}) per time-step (it), in order to obtain convergence within full machine accuracy, and the execution time in {\em sec}. From the listed results, one infers that:
\begin{itemize}
\item as is striking clear, the higher order methods are much more efficient than the lower order ones, especially when a high accuracy is required;

\item the theoretical rate of convergence for both the solution and the Hamiltonian errors is that we expected (for PHBVMs, until the Hamiltonian error falls within the round-off error level);

\item for a fixed time-step $h$, the numerical solutions provided by the Gauss methods and by the corresponding PHBVM method have a comparable accuracy, despite the negligible Hamiltonian error of these latter methods; 

\item the execution times of the PHBVM methods are about double than those of the corresponding Gauss methods, even though the mean number of blended iterations per time steps is practically the same (this latter, decreasing with the time-step $h$, and slightly increasing with $s$).

\end{itemize}
As a result, one would conclude that the conservation of the Hamiltonian apparently gains no practical advantage. However, this conclusion is readily confuted if we look at the error growth in the Hamiltonian and in the solution. In fact, in Figure~\ref{ex1_err} there is the plot of the Hamiltonian error (left plot) and the solution error (right plot) by using the
3-stage Gauss method and the PHBVM(6,3) method with time-step $h=T/100$ over 100 periods. As one may see, now it is clear that the 3-stage Gauss method exhibits a numerical drift in the energy, unlike PHBVM(6,3). As a result, this latter method exhibits a linear error growth, whereas the former one has a quadratic error growth.

\begin{table}[t]
\caption{results for problem (\ref{ex1}).}
\label{ex1tab}
\centering  
\scriptsize
\begin{tabular}{|r|rr|rr|r|r|rr|rr|r|r|}
\hline
 &\multicolumn{6}{|c|}{Gauss-1} &\multicolumn{6}{|c|}{PHBVM$(4,1)$}\\
\hline
$n$ & $e_y$ &rate& $e_H$ & rate & it & time& $e_y$ &rate& $e_H$ & rate & it & time\\
\hline
   50 & 3.54e-02 & --- & 4.47e-02 & --- &   7.4 &  0.01 & 7.64e-02 & --- & 1.72e-07 & --- &   8.5 &  0.02 \\  
  100 & 8.56e-03 & 2.1 & 1.09e-02 & 2.0 &   5.8 &  0.01 & 1.85e-02 & 2.0 & 6.48e-10 & 8.1 &   6.7 &  0.03 \\ 
  200 & 2.12e-03 & 2.0 & 2.71e-03 & 2.0 &   5.1 &  0.02 & 4.58e-03 & 2.0 & 2.37e-12 & 8.1 &   5.5 &  0.05 \\ 
  400 & 5.29e-04 & 2.0 & 6.77e-04 & 2.0 &   4.4 &  0.04 & 1.14e-03 & 2.0 & 8.88e-16 & ** &   4.6 &  0.09 \\ 
  800 & 1.32e-04 & 2.0 & 1.69e-04 & 2.0 &   4.1 &  0.07 & 2.86e-04 & 2.0 & 8.88e-16 & ** &   4.3 &  0.16 \\ 
 1600 & 3.30e-05 & 2.0 & 4.23e-05 & 2.0 &   3.5 &  0.13 & 7.14e-05 & 2.0 & 1.78e-15 & ** &   4.0 &  0.30 \\ 
 3200 & 8.25e-06 & 2.0 & 1.06e-05 & 2.0 &   3.2 &  0.25 & 1.79e-05 & 2.0 & 8.88e-16 & ** &   3.4 &  0.56 \\  
 6400 & 2.06e-06 & 2.0 & 2.64e-06 & 2.0 &   3.1 &  0.48 & 4.46e-06 & 2.0 & 1.78e-15 & ** &   3.1 &  1.08 \\
12800 & 5.16e-07 & 2.0 & 6.61e-07 & 2.0 &   3.0 &  0.95 & 1.12e-06 & 2.0 & 1.78e-15 & ** &   3.0 &  2.12 \\ 
25600 & 1.29e-07 & 2.0 & 1.65e-07 & 2.0 &   3.0 &  1.95 & 2.79e-07 & 2.0 & 1.78e-15 & ** &   3.0 &  4.18 \\ 
51200 & 3.22e-08 & 2.0 & 4.13e-08 & 2.0 &   3.0 &  3.87 & 6.97e-08 & 2.0 & 1.78e-15 & ** &   3.0 &  8.41 \\ 
102400 & 8.06e-09 & 2.0 & 1.03e-08 & 2.0 &   2.6 &  7.45 & 1.74e-08 & 2.0 & 1.78e-15 & ** &   2.9 & 16.54 \\ 
204800 & 2.02e-09 & 2.0 & 2.58e-09 & 2.0 &   2.3 & 14.32 & 4.36e-09 & 2.0 & 2.66e-15 & ** &   2.4 & 30.83 \\ 
409600 & 5.04e-10 & 2.0 & 6.45e-10 & 2.0 &   2.1 & 27.64 & 1.09e-09 & 2.0 & 2.66e-15 & ** &   2.2 & 59.74 \\  
819200 & 1.26e-10 & 2.0 & 1.61e-10 & 2.0 &   2.0 & 55.22 & 2.72e-10 & 2.0 & 2.66e-15 & ** &   2.0 & 114.18 \\ 
\hline
&\multicolumn{6}{|c|}{Gauss-2} &\multicolumn{6}{|c|}{PHBVM$(4,2)$}\\
\hline
$n$ & $e_y$ &rate& $e_H$ & rate & it & time& $e_y$ &rate& $e_H$ & rate & it & time\\
\hline
   50 & 3.43e-04 & --- & 1.83e-04 & --- &   8.9 &  0.01& 4.89e-05 & --- & 7.97e-09 & --- &   9.1 &  0.02 \\ 
  100 & 2.16e-05 & 4.0 & 1.15e-05 & 4.0 &   7.8 &  0.02& 3.05e-06 & 4.0 & 3.19e-11 & 8.0 &   7.9 &  0.03 \\ 
  200 & 1.35e-06 & 4.0 & 7.21e-07 & 4.0 &   6.8 &  0.03& 1.90e-07 & 4.0 & 8.88e-16 & ** &   6.9 &  0.05 \\ 
  400 & 8.44e-08 & 4.0 & 4.51e-08 & 4.0 &   6.0 &  0.06& 1.19e-08 & 4.0 & 8.88e-16 & ** &   6.1 &  0.10 \\ 
  800 & 5.28e-09 & 4.0 & 2.82e-09 & 4.0 &   5.5 &  0.10& 7.44e-10 & 4.0 & 1.78e-15 & ** &   5.6 &  0.18 \\ 
 1600 & 3.30e-10 & 4.0 & 1.76e-10 & 4.0 &   5.1 &  0.19& 4.65e-11 & 4.0 & 8.88e-16 & ** &   5.2 &  0.35 \\ 
 3200 & 2.06e-11 & 4.0 & 1.10e-11 & 4.0 &   4.7 &  0.36& 2.97e-12 & 4.0 & 1.78e-15 & ** &   4.7 &  0.65 \\ 
 6400 & 1.28e-12 & 4.0 & 6.84e-13 & 4.0 &   4.3 &  0.68& 2.24e-13 & 3.7 & 1.78e-15 & ** &   4.3 &  1.23 \\ 
\hline
&\multicolumn{6}{|c|}{Gauss-3} &\multicolumn{6}{|c|}{PHBVM$(6,3)$}\\
\hline
$n$ & $e_y$ &rate& $e_H$ & rate & it & time& $e_y$ &rate& $e_H$ & rate & it & time\\
\hline
   50 & 5.49e-07 & --- & 2.88e-07 & --- &   9.7 &  0.01& 1.23e-07 & --- & 8.88e-16 & --- &   9.8 &  0.02 \\ 
  100 & 8.58e-09 & 6.0 & 4.49e-09 & 6.0 &   8.1 &  0.02& 1.92e-09 & 6.0 & 8.88e-16 &** &   8.2 &  0.03 \\ 
  200 & 1.34e-10 & 6.0 & 7.00e-11 & 6.0 &   7.0 &  0.03& 3.00e-11 & 6.0 & 8.88e-16 & ** &   7.1 &  0.06 \\ 
  400 & 2.12e-12 & 6.0 & 1.10e-12 & 6.0 &   6.3 &  0.06& 5.08e-13 & 5.9 & 1.78e-15 & ** &   6.4 &  0.12 \\ 
  800 & 5.30e-14 & 5.3 & 2.04e-14 & 5.7 &   5.7 &  0.11& 4.80e-14 & 3.4 & 1.78e-15 & ** &   5.7 &  0.22 \\ 
\hline
\end{tabular}
\end{table}

\begin{figure}[t]
\centerline{\includegraphics[width=9.5cm]{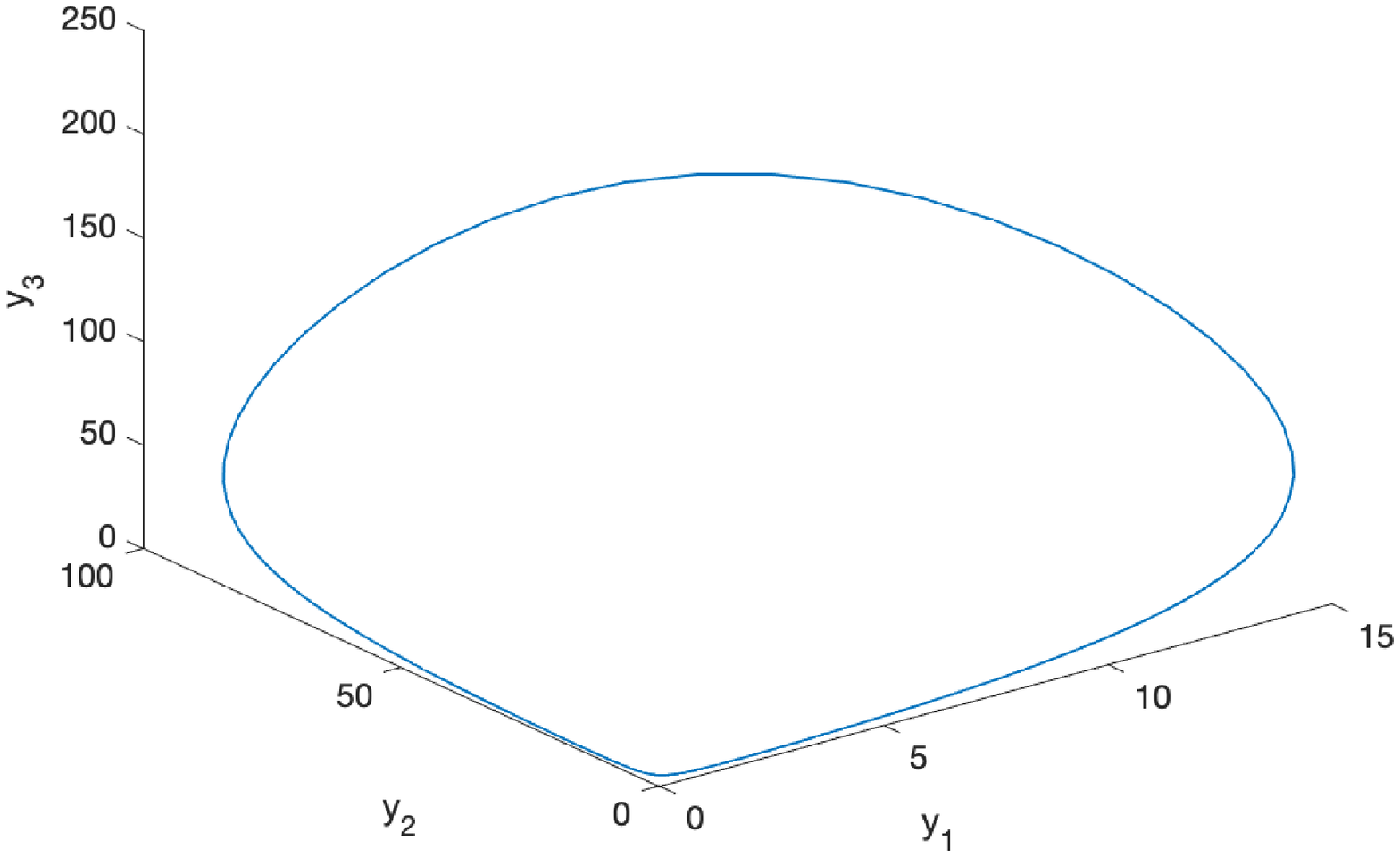}}
\caption{solution of problem (\ref{ex3}).}
\label{fig3}
\end{figure}

\bigskip
{\bf Exampe~2.~} The second example, taken from \cite{GNI2006}, is given by:   

\begin{eqnarray}\nonumber
\dot y &=& \pmatrix{ccc} 0 & y_1y_2 & y_1y_3\\ -y_1y_2 &0 & -y_2y_3\\ -y_1y_3 & y_2y_3 & 0\endpmatrix \nabla H(y),\\[2mm] \label{ex3}
H(y)   &=& a\left(\ln y_1 -\frac{y_1}{y_1^*}\right) +b\left(\ln y_2 -\frac{y_2}{y_2^*}\right) + c\left(\ln y_3 -\frac{y_3}{y_3^*}\right),\\[2mm]
C(y)  &=& -\ln y_1 -\ln y_2 +\ln y_3,\nonumber
\end{eqnarray}
with
$$a=1, \qquad b=2, \qquad c=3, \qquad y_1^*=1, \qquad y_2^*=10,\qquad y_3^* = 50,\qquad y(0) = \left(1,\,1,\,1\right)^\top,$$
whose solution, which is periodic of period $T\approx 2.143610709155912$, is depicted in Figure~\ref{fig3}.

\begin{figure}[t]
\centerline{\includegraphics[width=7.25cm]{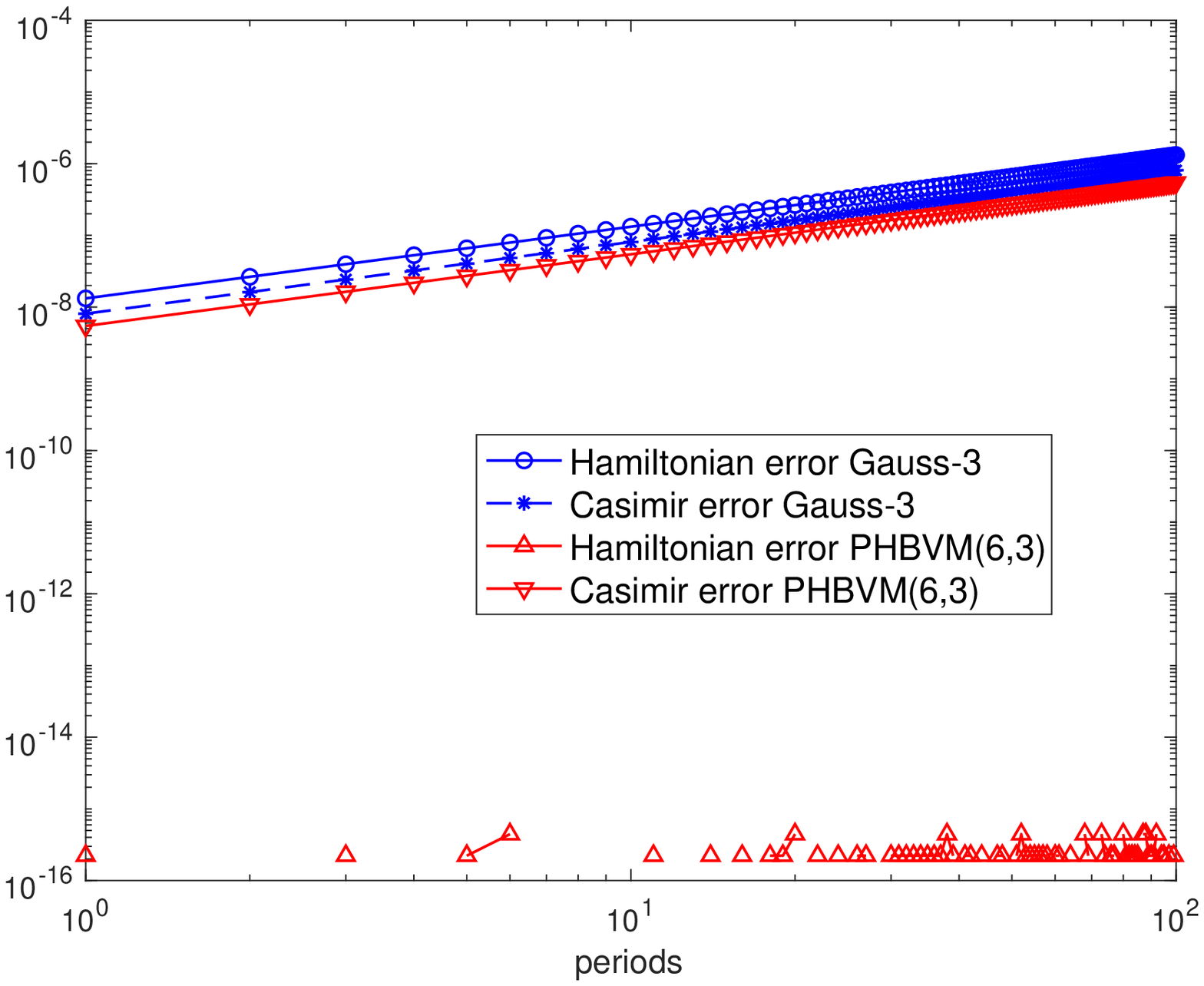}\qquad\includegraphics[width=7.25cm]{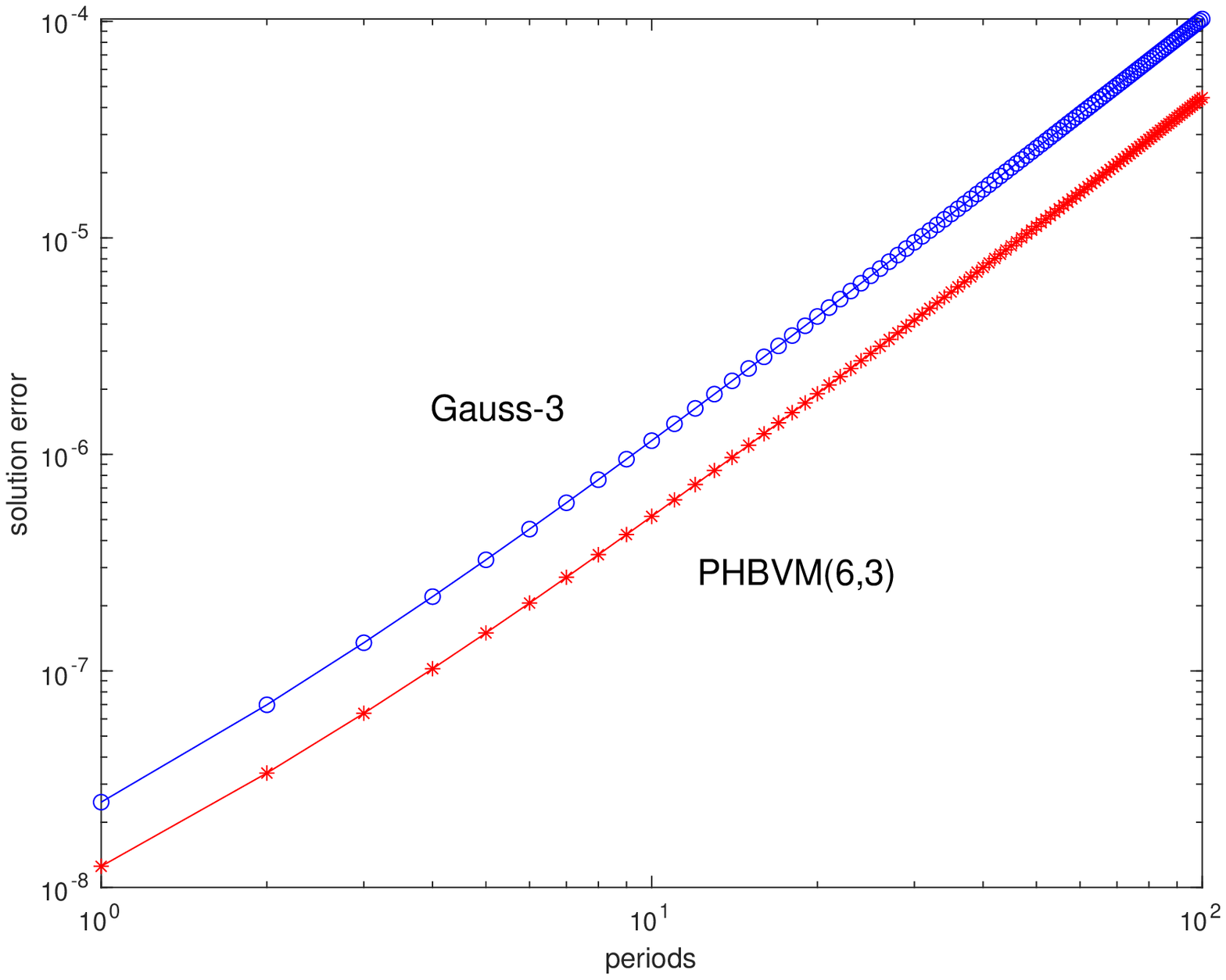}}
\caption{Hamiltonian and Casimir errors (left plot) and solution error (right plot) when solving problem (\ref{ex3}) with time-step $h=T/100$ over 100 periods with the 3-stage Gauss and PHBVM(6,3) methods.}
\label{ex3_err}
\end{figure}

\begin{figure}[t]
\centerline{\includegraphics[width=9.5cm]{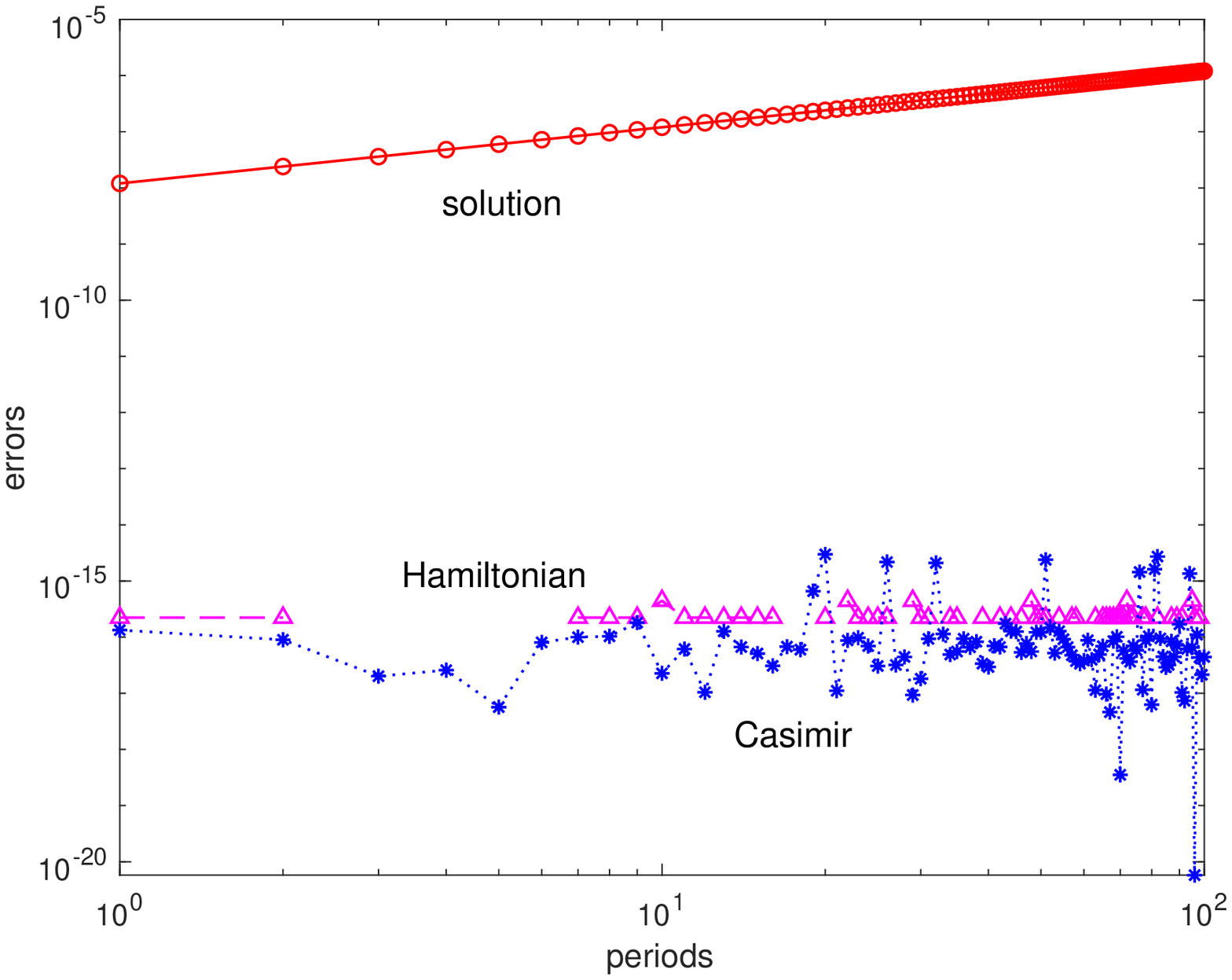}}
\caption{Hamiltonian, Casimir, and solution errors when solving problem (\ref{ex3}) with time-step $h=T/100$ over 100 periods with the EPHBVM(6,3) method.}
\label{ex3_err1}
\end{figure}

At first, we compare the same methods used for the previous example, again with time-step $h=T/n$. The obtained results for the $s$-stage Gauss and PHBVM methods
are listed in Table~\ref{ex3tab}: the conclusions that one can derive from them are similar to those driven from Table~\ref{ex1tab} for the previous example, with the additional remark that now the Casimir $C(y)$ is not conserved.\footnote{In Table~\ref{ex3tab} $e_C$ denotes the error in the Casimir.} This fact, in turn, produces the results depicted in the two plots in Figure~\ref{ex3_err}, concerning the application of the Gauss-3 and PHBVM(6,3) methods for solving the problem with time-step $h=T/100$ over 100 periods. From the two plots, one infers that both methods exhibit a drift in the Casimir, whereas only Gauss-3 exhibits a drift in the Hamiltonian, too (left plot); however, {\em both} methods exhibit a quadratic error growth in the solution, despite the fact that PHBVM(6,3) conserves the Hamiltonian.
For this reason, in Table~\ref{ex3tab1} we list the obtained results by using the EHBVM(4,1), EHBVM(4,2), and EHBVM(6,3) methods for solving problem (\ref{ex3}), by using the same time-steps considered for obtaining the results of Table~\ref{ex3tab}. As one may see, now the conservation of the Casimir is soon obtained, as the time-step is decreased, besides that of the Hamiltonian, with a computational cost perfectly comparable to that of the corresponding PHBVM method. The conservation of both invariants, in turn, allows to recover a linear error growth in the numerical solution, as is shown in the plot of Figure~\ref{ex3_err1}.

\begin{sidewaystable}
\caption{results for problem (\ref{ex3}).}
\label{ex3tab}
\centering
\scriptsize
\begin{tabular}{|r|rr|rr|rr|r|r|rr|rr|rr|r|r|}
\hline
 &\multicolumn{8}{|c|}{Gauss-1} &\multicolumn{8}{|c|}{PHBVM$(4,1)$}\\
\hline
$n$ & $e_y$ &rate& $e_H$ & rate & $e_C$ & rate & it & time& $e_y$ &rate& $e_H$ & rate & $e_C$ & rate& it & time\\
\hline
    50 & 6.25e-02 & ---& 4.89e-01 & ---& 3.86e-02 & ---&   8.4 &  0.03& 1.23e-01 & ---& 1.01e-05 & ---& 5.45e-02 & ---&   9.8 &  0.05 \\ 
  100 & 1.62e-02 & 1.9 & 1.26e-01 & 2.0 & 9.80e-03 & 2.0 &   6.5 &  0.01& 3.00e-02 & 2.0 & 3.80e-08 & 8.0 & 1.32e-02 & 2.0 &   7.2 &  0.03 \\ 
  200 & 4.10e-03 & 2.0 & 3.18e-02 & 2.0 & 2.44e-03 & 2.0 &   5.6 &  0.03& 7.46e-03 & 2.0 & 5.55e-15 & ** & 3.26e-03 & 2.0 &   5.9 &  0.06 \\ 
  400 & 1.03e-03 & 2.0 & 7.97e-03 & 2.0 & 6.09e-04 & 2.0 &   4.7 &  0.04& 1.86e-03 & 2.0 & 5.11e-15 & ** & 8.13e-04 & 2.0 &   5.0 &  0.10 \\ 
  800 & 2.57e-04 & 2.0 & 1.99e-03 & 2.0 & 1.52e-04 & 2.0 &   4.3 &  0.08& 4.65e-04 & 2.0 & 2.00e-15 & ** & 2.03e-04 & 2.0 &   4.4 &  0.18 \\ 
 1600 & 6.42e-05 & 2.0 & 4.99e-04 & 2.0 & 3.81e-05 & 2.0 &   3.8 &  0.14& 1.16e-04 & 2.0 & 2.44e-15 & ** & 5.07e-05 & 2.0 &   4.1 &  0.33 \\ 
 3200 & 1.60e-05 & 2.0 & 1.25e-04 & 2.0 & 9.51e-06 & 2.0 &   3.5 &  0.25& 2.91e-05 & 2.0 & 3.77e-15 & ** & 1.27e-05 & 2.0 &   3.6 &  0.59 \\ 
 6400 & 4.01e-06 & 2.0 & 3.12e-05 & 2.0 & 2.38e-06 & 2.0 &   3.1 &  0.49& 7.27e-06 & 2.0 & 3.33e-15 & ** & 3.17e-06 & 2.0 &   3.4 &  1.14 \\ 
12800 & 1.00e-06 & 2.0 & 7.79e-06 & 2.0 & 5.95e-07 & 2.0 &   3.0 &  0.97& 1.82e-06 & 2.0 & 3.77e-15 & ** & 7.93e-07 & 2.0 &   3.1 &  2.19 \\ 
25600 & 2.51e-07 & 2.0 & 1.95e-06 & 2.0 & 1.49e-07 & 2.0 &   3.0 &  1.93& 4.54e-07 & 2.0 & 3.77e-15 & ** & 1.98e-07 & 2.0 &   3.0 &  4.32 \\ 
51200 & 6.27e-08 & 2.0 & 4.87e-07 & 2.0 & 3.72e-08 & 2.0 &   3.0 &  3.87& 1.14e-07 & 2.0 & 3.77e-15 & ** & 4.95e-08 & 2.0 &   3.0 &  8.61 \\ 
102400 & 1.57e-08 & 2.0 & 1.22e-07 & 2.0 & 9.29e-09 & 2.0 &   2.9 &  7.67& 2.84e-08 & 2.0 & 4.22e-15 & ** & 1.24e-08 & 2.0 &   3.0 & 17.18 \\ 
204800 & 3.92e-09 & 2.0 & 3.04e-08 & 2.0 & 2.32e-09 & 2.0 &   2.5 & 14.70& 7.10e-09 & 2.0 & 5.11e-15 & ** & 3.10e-09 & 2.0 &   2.7 & 33.34 \\ 
409600 & 9.79e-10 & 2.0 & 7.61e-09 & 2.0 & 5.81e-10 & 2.0 &   2.3 & 28.79& 1.77e-09 & 2.0 & 5.55e-15 & ** & 7.74e-10 & 2.0 &   2.5 & 63.85 \\ 
819200 & 2.45e-10 & 2.0 & 1.90e-09 & 2.0 & 1.45e-10 & 2.0 &   2.1 & 55.80& 4.44e-10 & 2.0 & 5.55e-15 & ** & 1.93e-10 & 2.0 &   2.2 & 122.89 \\ 
\hline
&\multicolumn{8}{|c|}{Gauss-2} &\multicolumn{8}{|c|}{PHBVM$(4,2)$}\\
\hline
$n$ & $e_y$ &rate& $e_H$ & rate & $e_C$ & rate & it & time& $e_y$ &rate& $e_H$ & rate & $e_C$ & rate& it & time\\
\hline
   50 & 2.56e-04 & ---& 1.79e-03 & ---& 8.08e-04 & ---&  10.4 &  0.02& 2.18e-04 & ---& 3.49e-07 & ---& 9.72e-04 & ---&  10.4 &  0.03 \\ 
  100 & 1.58e-05 & 4.0 & 1.11e-04 & 4.0 & 5.37e-05 & 3.9 &   8.5 &  0.02& 1.30e-05 & 4.1 & 1.52e-09 & 7.8 & 6.22e-05 & 4.0 &   8.7 &  0.04 \\ 
  200 & 9.88e-07 & 4.0 & 6.97e-06 & 4.0 & 3.32e-06 & 4.0 &   7.4 &  0.03& 8.05e-07 & 4.0 & 3.77e-15 & ** & 3.86e-06 & 4.0 &   7.5 &  0.06 \\ 
  400 & 6.17e-08 & 4.0 & 4.35e-07 & 4.0 & 2.09e-07 & 4.0 &   6.5 &  0.06& 5.02e-08 & 4.0 & 2.00e-15 & ** & 2.42e-07 & 4.0 &   6.5 &  0.11 \\ 
  800 & 3.86e-09 & 4.0 & 2.72e-08 & 4.0 & 1.31e-08 & 4.0 &   6.0 &  0.11& 3.13e-09 & 4.0 & 2.00e-15 & ** & 1.51e-08 & 4.0 &   6.2 &  0.21 \\ 
 1600 & 2.41e-10 & 4.0 & 1.70e-09 & 4.0 & 8.16e-10 & 4.0 &   5.4 &  0.20& 1.96e-10 & 4.0 & 3.77e-15 & ** & 9.46e-10 & 4.0 &   5.4 &  0.37 \\ 
 3200 & 1.51e-11 & 4.0 & 1.06e-10 & 4.0 & 5.10e-11 & 4.0 &   5.0 &  0.38& 1.23e-11 & 4.0 & 3.77e-15 & ** & 5.91e-11 & 4.0 &   5.1 &  0.71 \\ 
 6400 & 9.98e-13 & 3.9 & 6.65e-12 & 4.0 & 3.19e-12 & 4.0 &   4.5 &  0.71& 8.54e-13 & 3.8 & 3.77e-15 & ** & 3.70e-12 & 4.0 &   4.5 &  1.35 \\ 
\hline
&\multicolumn{8}{|c|}{Gauss-3} &\multicolumn{8}{|c|}{PHBVM$(6,3)$}\\
\hline
$n$ & $e_y$ &rate& $e_H$ & rate & $e_C$ & rate & it & time& $e_y$ &rate& $e_H$ & rate & $e_C$ & rate& it & time\\
\hline
   50 & 1.13e-06 & --- & 6.47e-06 & --- & 3.57e-06 & --- &  10.9 &  0.01& 5.51e-07 & --- & 5.11e-15 & --- & 1.97e-06 & --- &  11.1 &  0.03 \\ 
  100 & 1.76e-08 & 6.0 & 1.00e-07 & 6.0 & 5.43e-08 & 6.0 &   8.9 &  0.02& 9.34e-09 & 5.9 & 3.33e-15 & ** & 2.79e-08 & 6.1 &   9.1 &  0.04 \\ 
  200 & 2.76e-10 & 6.0 & 1.60e-09 & 6.0 & 8.86e-10 & 5.9 &   7.7 &  0.03& 1.49e-10 & 6.0 & 2.89e-15 & ** & 4.26e-10 & 6.0 &   7.8 &  0.07 \\ 
  400 & 4.23e-12 & 6.0 & 2.50e-11 & 6.0 & 1.38e-11 & 6.0 &   6.8 &  0.06& 2.25e-12 & 6.0 & 2.89e-15 & ** & 6.62e-12 & 6.0 &   6.9 &  0.12 \\ 
  800 & 3.08e-14 & 7.1 & 4.05e-13 & 5.9 & 2.17e-13 & 6.0 &   6.3 &  0.11& 6.30e-14 & 5.2 & 3.33e-15 & ** & 1.03e-13 & 6.0 &   6.3 &  0.23 \\ 
\hline
\end{tabular}
\end{sidewaystable}

\begin{table}[t]
\caption{further results for problem (\ref{ex3}).}
\label{ex3tab1}
\centering
\scriptsize
\begin{tabular}{|r|rr|rr|rr|r|r|}
\hline
 &\multicolumn{8}{|c|}{EPHBVM$(4,1)$}\\
\hline
$n$ & $e_y$ &rate& $e_H$ & rate & $e_C$ & rate & it & time\\
\hline
   50 & 1.26e-01 & --- & 9.36e-06 & --- & 2.21e-06 & --- &  11.1 &  0.06 \\ 
  100 & 3.07e-02 & 2.0 & 3.73e-08 & 8.0 & 9.48e-09 & 7.9 &   8.2 &  0.04 \\ 
  200 & 7.45e-03 & 2.0 & 5.55e-15 & ** & 1.78e-15 & ** &   6.8 &  0.07 \\ 
  400 & 1.90e-03 & 2.0 & 3.77e-15 & ** & 1.78e-15 & ** &   5.6 &  0.12 \\ 
  800 & 4.70e-04 & 2.0 & 3.77e-15 & ** & 8.88e-16 & ** &   4.8 &  0.22 \\ 
 1600 & 1.18e-04 & 2.0 & 3.77e-15 & ** & 1.78e-15 & ** &   4.4 &  0.39 \\ 
 3200 & 2.96e-05 & 2.0 & 3.77e-15 & ** & 1.78e-15 & ** &   4.1 &  0.71 \\ 
 6400 & 7.44e-06 & 2.0 & 5.11e-15 & ** & 1.78e-15 & ** &   3.9 &  1.39 \\ 
12800 & 1.83e-06 & 2.0 & 5.11e-15 & ** & 1.78e-15 & ** &   3.5 &  2.64 \\ 
25600 & 4.64e-07 & 2.0 & 6.88e-15 & ** & 1.78e-15 & ** &   3.3 &  5.16 \\ 
51200 & 1.18e-07 & 2.0 & 6.88e-15 & ** & 1.78e-15 & ** &   3.1 &  9.98 \\ 
102400 & 2.83e-08 & 2.1 & 6.88e-15 & ** & 1.78e-15 & ** &   3.0 & 19.71 \\ 
204800 & 7.20e-09 & 2.0 & 5.55e-15 & ** & 1.78e-15 & ** &   2.8 & 38.16 \\ 
409600 & 1.79e-09 & 2.0 & 6.88e-15 & ** & 1.78e-15 & ** &   2.5 & 72.98 \\ 
819200 & 4.53e-10 & 2.0 & 7.33e-15 & ** & 1.78e-15 & ** &   2.2 & 140.73 \\
\hline
&\multicolumn{8}{|c|}{EPHBVM$(4,2)$}\\
\hline
$n$ & $e_y$ &rate& $e_H$ & rate & $e_C$ & rate & it & time\\
\hline
   50 & 1.29e-04 & --- & 3.49e-07 & --- & 1.72e-07 & --- &  10.4 &  0.03 \\ 
  100 & 6.20e-06 & 4.4 & 1.52e-09 & 7.8 & 6.01e-10 & 8.2 &   8.7 &  0.04 \\ 
  200 & 1.34e-07 & 5.5 & 5.11e-15 & ** & 8.88e-16 & ** &   7.5 &  0.07 \\ 
  400 & 2.54e-08 & 2.4 & 3.33e-15 & ** & 1.78e-15 & ** &   6.5 &  0.12 \\ 
  800 & 6.73e-10 & 5.2 & 3.33e-15 & ** & 8.88e-16 & ** &   6.2 &  0.24 \\ 
 1600 & 9.47e-11 & 2.8 & 3.33e-15 & ** & 1.78e-15 & ** &   5.4 &  0.42 \\ 
 3200 & 6.25e-12 & 3.9 & 4.22e-15 & ** & 1.78e-15 & ** &   5.1 &  0.83 \\ 
 6400 & 5.46e-13 & 3.5 & 3.77e-15 & ** & 1.78e-15 & ** &   4.5 &  1.54 \\ 
\hline
&\multicolumn{8}{|c|}{EPHBVM$(6,3)$}\\
\hline
$n$ & $e_y$ &rate& $e_H$ & rate & $e_C$ & rate & it & time\\
\hline
   50 & 5.32e-07 & -- & 3.33e-15 & --- & 8.88e-16 & --- &  11.1 &  0.03 \\ 
  100 & 7.44e-09 & 6.2 & 3.77e-15 & ** & 1.78e-15 & ** &   9.1 &  0.05 \\ 
  200 & 5.78e-11 & 7.0 & 3.33e-15 & ** & 1.78e-15 & ** &   7.8 &  0.08 \\ 
  400 & 1.81e-12 & 5.0 & 2.00e-15 & ** & 8.88e-16 & ** &   6.9 &  0.14 \\ 
  800 & 8.24e-14 & 4.5 & 4.22e-15 & ** & 1.33e-15 & ** &   6.3 &  0.27 \\ 
\hline
\end{tabular}
\end{table}

\section{Conclusions}\label{fine}

In this paper we have presented a class of energy-conserving line integral methods for Poisson problems. In the case where the problem is Hamiltonian, these methods reduce to the class of Hamiltonian Boundary Value Methods (HBVMs), which are energy-conserving methods for such problems. Consequently, the new methods  can be regarded as an extension of HBVMs for Poisson (not Hamiltonian) problems, which we called PHBVMs. Moreover, a further enhancement of such methods (EPHBVMs) allows to obtain the conservation of Casimirs, too. A thorough analysis of the methods has been carried out,  confirmed by a couple of numerical tests. {\red As a further direction of investigation, we mention the study of the application of the  methods for solving highly-oscillatory Poisson problems, similarly as done with HBVMs in the Hamiltonian case \cite{osc1,osc2,osc3}.} 





\begin{thebibliography}{10}

\bibitem{osc3} P.\,Amodio, L.\,Brugnano, F.\,Iavernaro. Analysis of Spectral Hamiltonian Boundary Value Methods (SHBVMs) for the numerical solution of ODE problems.  {\em Numer. Algorithms}   {\bf 83} (2020) 1489--1508.

\bibitem{BCMR2012} L.\,Brugnano, M.\,Calvo, J.I.\,Montijano, L.\,R\'andez. Energy preserving methods for Poisson systems. {\em J. Comput. Appl. Math.} {\bf 236} (2012) 3890--3904.

\bibitem{BFCI2014}  L.\,Brugnano, G.\,Frasca Caccia, F.\,Iavernaro. Efficient  implementation of Gauss collocation and Hamiltonian Boundary Value Methods. {\em Numer. Algorithms} {\bf 65} (2014) 633--650.

\bibitem{BGI2018} L.\,Brugnano, G.\,Gurioli, F.\,Iavernaro. Analysis of Energy and QUadratic Invariant Preserving (EQUIP) methods.  {\em J. Comput. Appl. Math.} {\bf 335} (2018) 51--73.

\bibitem{BI2012} 
L.\,Brugnano, F.\,Iavernaro. Line integral methods which preserve all invariants of conservative problems. {\em J. Comput. Appl. Math.} {\bf 236} (2012) 3905--3919..

\bibitem{LIMbook2016} L.\,Brugnano, F.\,Iavernaro. {\em Line Integral Methods for Conservative Problems}.  Chapman and Hall/CRC, Boca Raton, FL, 2016.

\bibitem{BI2018} L.\,Brugnano, F.\,Iavernaro. Line Integral Solution of Differential Problems. {\em Axioms} {\bf 7}(2) (2018) article n.\,36. \url{http://dx.doi.org//10.3390/axioms7020036}

\bibitem{BIT2009} L.\,Brugnano, F.\,Iavernaro, D.\,Trigiante. Hamiltonian BVMs (HBVMs): A family of ``drift-free'' methods for integrating polynomial Hamiltonian systems. {\em AIP Conf. Proc.} {\bf 1168} (2009) 715--718.

\bibitem{BIT2010} L.\,Brugnano, F.\,Iavernaro, D.\,Trigiante.  Hamiltonian Boundary Value Methods (Energy Preserving Discrete Line Integral Methods).  {\em JNAIAM J. Numer. Anal. Ind. Appl. Math.} {\bf 5},\,1-2 (2010) 17--37.

\bibitem{BIT2011} L.\,Brugnano, F.\,Iavernaro, D.\,Trigiante. A note on the efficient implementation of Hamiltonian BVMs. {\em J. Comput. Appl. Math.} {\bf 236} (2011) 375--383.

\bibitem{BIT2012-1} L.\,Brugnano, F.\,Iavernaro, D.\,Trigiante.  The lack of continuity and the role of infinite and infinitesimal in numerical methods for ODEs: the case of symplecticity. {\em Appl. Math. Comput.} {\bf 218} (2012) 8056--8063.

\bibitem{BIT2012} L.\,Brugnano, F.\,Iavernaro, D.\,Trigiante.  A simple framework for the derivation and analysis of effective one-step methods for ODEs. {\em Appl. Math. Comput.} {\bf 218} (2012) 8475--8485.

\bibitem{BIT2015}  L.\,Brugnano, F.\,Iavernaro, D.\,Trigiante. Analisys of Hamiltonian Boundary Value Methods (HBVMs): A class of energy-preserving Runge-Kutta methods for the numerical solution of polynomial Hamiltonian systems. {\em Commun. Nonlinear Sci. Numer. Simul.} {\bf 20} (2015) 650--667. 

\bibitem{osc2} L.\,Brugnano, F.\,Iavernaro, J.I.\,Montijano, L.\,R\'andez. Spectrally accurate space-time solution of Hamiltonian PDEs. {\em Numer. Algorithms}   {\bf 81} (2019) 1183--1202.

\bibitem{BM2002} L.\,Brugnano, C.\,Magherini. Blended Implementation of Block Implicit Methods for ODEs. {\em Appl. Numer. Math.} {\bf 42} (2002) 29--45.

\bibitem{BM2004} L.\,Brugnano, C.\,Magherini. The BiM Code for the Numerical Solution of ODEs.  {\em J. Comput. Appl. Math.} {\bf 164-165} (2002) 145--158.

\bibitem{BMM2006} L.\,Brugnano, C.\,Magherini.  Blended Implicit Methods for the Numerical Solution of DAE Problems.  {\em J. Comput. Appl. Math.} {\bf 189} (2006) 34--50.

\bibitem{BM2009} L.\,Brugnano, C.\,Magherini. Recent advances in linear analysis of convergence for splittings for solving ODE problems. {\em Appl. Numer. Math.} {\bf 59} (2009) 542--557.

\bibitem{osc1} L.\,Brugnano, J.I.\,Montijano, L.\,R\'andez. On the effectiveness of spectral methods for the numerical solution of multi-frequency highly-oscillatory Hamiltonian problems.  {\em Numer. Algorithms}   {\bf 81} (2019) 345--376.

\bibitem{BMR2019}  L.\,Brugnano, J.I.\,Montijano, L.\,R\'andez. High-order energy-conserving Line Integral Methods for charged particle dynamics.  {\em J. Comput. Phys.} {\bf 396} (2019) 209--227.

\bibitem{BS2014} L.\,Brugnano, Y.\,Sun. Multiple invariants conserving Runge-Kutta type methods for Hamiltonian problems. {\em Numer. Algorithms} {\bf 65} (2014) 611--632.

\bibitem{CH2011} D.\,Cohen, E.\,Hairer. Linear energy-preserving integrators for Poisson systems. {\em BIT Numer. Math.} {\bf 51} (2011) 91--101.


\bibitem{GNI2006} E.\,Hairer, C.\,Lubich, G.\,Wanner. {\em Geometric Numerical Integration, 2nd ed.} Springer, Berlin, 2006.

\bibitem{MHW2022} L.\,Mei, L.\,Huang, X.\,Wu. A unified framework for the study of high-order energy-preserving integrators for solving Poisson systems. {\em J. Comput. Phys.} \url{https://doi.org/10.1016/j.jcp.2021.110822}

\bibitem{M2015} Y.\,Miyatake. A derivation of energy-preserving exponentially-fitted integrators for Poisson systems. {\em  Comput. Phys. Commun.} {\bf 187} (2015) 156--161.



\bibitem{QML2008} G.R.W.\,Quispel,  D.I.\,McLaren. A new class of energy-preserving numerical integration methods. {\em J. Phys. A} {\bf 41}(4) (2008) 045206.

\bibitem{WMF2017} B.\,Wang, F.\,Meng, Y.\,Fang. Efficient implementation of RKN-type Fourier collocation methods for second-order differential equations.  {\em Appl. Numer. Math.} {\bf 119} (2017) 164--178. 

\bibitem{WW2018} B.\,Wang, X.\,Wu. Functionally-fitted energy-preserving integrators for Poisson systems. {\em J. Comput. Phys.} {\bf 364} (2018) 137--152.

\red
\bibitem{WW2021} B.\,Wang, X.\,Wu. {\em Geometric Integrators for Differential Equations with Highly Oscillatory Solutions.} Springer Nature Singapore Pte Ltd., 2021.

\end{thebibliography}
\end{document}